\documentclass[12pt]{article}
\usepackage{amsmath,amssymb}
\newcommand{\beq}{\begin{eqnarray}}
\newcommand{\eeq}{\end{eqnarray}}
\newtheorem{theorem}{Theorem}

\newtheorem{lemma}{Lemma}

\newcommand{\proof}{\medskip\par\noindent {\it Proof.\/}\quad}
\newcommand{\qed}{{\it Q.E.D.\/} \bigskip\par}
\newcommand{\remark}{\medskip\par\noindent {\it Remark.\/}\quad}
\newcommand{\rd}{\partial}
\newcommand{\const}{\mathrm{constant}}
\renewcommand{\Im}{\mathrm{Im}\,}

\newcommand{\bbC}{\mathbb{C}}
\newcommand{\bbP}{\mathbb{P}}

\newcommand{\bbZ}{\mathbb{Z}}
\newcommand{\calH}{\mathcal{H}}
\newcommand{\calL}{\mathcal{L}}
\newcommand{\PI}{\mathrm{P_I}}
\newcommand{\PII}{\mathrm{P_{II}}}
\newcommand{\PIII}{\mathrm{P_{III}}}
\newcommand{\PIV}{\mathrm{P_{IV}}}
\newcommand{\PV}{\mathrm{P_V}}
\newcommand{\PVI}{\mathrm{P_{VI}}}
\begin{document}

%%%%%%%%%%%%%%%%%%%%
%%%% title page %%%%
%%%%%%%%%%%%%%%%%%%%
\title{Painlev\'e-Calogero Correspondence Revisited}
\author{Kanehisa Takasaki\\
{\normalsize Department of Fundamental Sciences, Kyoto University}\\
{\normalsize Yoshida, Sakyo-ku, Kyoto 606-8501, Japan}\\
{\normalsize E-mail: takasaki@math.h.kyoto-u.ac.jp}}
\date{}
\maketitle

\begin{abstract}
We extend the work of Fuchs, Painlev\'e and Manin on 
a Calogero-like expression of the sixth Painlev\'e 
equation (the ``Painlev\'e-Calogero correspondence'') 
to the other five Painlev\'e equations.  The Calogero 
side of the sixth Painlev\'e equation is known to be 
a non-autonomous version of the (rank one) elliptic 
model of Inozemtsev's extended Calogero systems.  
The fifth and fourth Painlev\'e equations correspond 
to the hyperbolic and rational  models in Inozemtsev's 
classification.  Those corresponding to the third, 
second and first are not included therein.  We further 
extend the correspondence to the higher rank models, 
and obtain a ``multi-component'' version of the 
Painlev\'e equations.  
\end{abstract}
\begin{flushleft}
\end{flushleft}
\bigskip

%%%%%%%%%%%%%%%%%%%%%%%%%%
%%% preprint data %%%%%%%%
%%%%%%%%%%%%%%%%%%%%%%%%%%
\begin{flushleft}
KUCP 149\\
{\tt math.QA/}0004118
\end{flushleft}

%%%%%%%%%%%%%%%%%%%
%%%% main text %%%%
%%%%%%%%%%%%%%%%%%%
\newpage
\section{Introduction}

The so called Painlev\'e equations are 
the following six equations discovered 
by Painlev\'e \cite{bib:Panleve} and 
Gambier \cite{bib:Gambier}:
\beq
  &(\PVI)& 
  \frac{d^2\lambda}{dt^2} 
  = \frac{1}{2}\left(\frac{1}{\lambda} 
    + \frac{1}{\lambda-1} + \frac{1}{\lambda-t}\right) 
    \left(\frac{d\lambda}{dt}\right)^2 
  - \left(\frac{1}{t} + \frac{1}{t-1} 
    + \frac{1}{\lambda-t}\right)\frac{d\lambda}{dt} 
   \nonumber \\
        &&
  + \frac{\lambda(\lambda-1)(\lambda-t)}{t^2(t-1)^2}
    \left(\alpha + \frac{\beta t}{\lambda^2} 
    + \frac{\gamma(t-1)}{(\lambda-1)^2} 
    + \frac{\delta t(t-1)}{(\lambda-t)^2}\right). 
   \nonumber \\
  &(\PV)& 
  \frac{d^2\lambda}{dt^2} 
  = \left(\frac{1}{2\lambda} + \frac{1}{\lambda-1}
    \right)\left(\frac{d\lambda}{dt}\right)^2 
  - \frac{1}{t}\frac{d\lambda}{dt} 
   \nonumber \\
       && 
  + \frac{\lambda(\lambda-1)^2}{t^2} 
    \left(\alpha + \frac{\beta}{\lambda^2} 
    + \frac{\gamma t}{(\lambda-1)^2} 
    + \frac{\delta t^2(\lambda+1)}{(\lambda-1)^3}\right). 
   \nonumber \\
  &(\PIV)&  
  \frac{d^2\lambda}{dt^2} 
  = \frac{1}{2\lambda}\left(\frac{d\lambda}{dt}\right)^2 
  + \frac{3}{2}\lambda^3 + 4t\lambda^2 
  + 2(t^2 - \alpha)\lambda + \frac{\beta}{\lambda}. 
    \nonumber \\
  &(\PIII)& 
  \frac{d^2\lambda}{dt^2} 
  = \frac{1}{\lambda}\left(\frac{d\lambda}{dt}\right)^2 
  - \frac{1}{t}\frac{d\lambda}{dt} 
  + \frac{\lambda^2}{4t^2}\left(\alpha 
    + \frac{\beta t}{\lambda^2} + \gamma\lambda 
    + \frac{\delta t^2}{4\lambda^3}\right). 
   \nonumber \\
  &(\PII)& 
  \frac{d^2\lambda}{dt^2} 
  = 2\lambda^3 + t\lambda + \alpha. 
   \nonumber \\ 
  &(\PI)& 
  \frac{d^2\lambda}{dt^2} 
  = 6\lambda^2 + t. 
\nonumber 
\eeq
The third equation $\PIII$ is slightly modified; 
the original equation can be reproduced by 
the simple change of variables $(t,\lambda) 
\to (t^2, t\lambda)$.  It is well known that 
these equations are characterized by the 
absence of ``movable singularities'' other than 
poles.  

R. Fuchs \cite{bib:Fuchs} proposed two more approaches 
to the sixth equation $\PVI$.  One approach is the 
concept of isomonodromic deformations. In this approach, 
$\PVI$ is interpreted as a differential equation 
describing isomonodromic deformations of a linear 
ordinary differential equation on the Riemann sphere.  
This is the origin of many subsequent researches.  
Another approach relates $\PVI$ to an incomplete 
elliptic integral.   Painlev\'e \cite{bib:Pa-elliptic} 
took the second approach, and derived a new expression 
of $\PVI$ in term of the Weierstrass $\wp$-function.   
This work of Painlev\'e is 
briefly reviewed in Okamoto's work on affine Weyl 
group symmetries of $\PVI$ \cite{bib:Ok-symmetry}. 

Manin \cite{bib:Manin} revived the almost 
forgotten work of Fuchs and Painlev\'e after 
nearly ninety years.  Manin's remarkable idea 
is to use the elliptic modulus $\tau$, rather 
than $t$, as an independent variable.  The 
outcome is a Hamiltonian system with a 
Hamiltonian of the normal form 
$\calH = p^2/2 + V(q)$, where the potential 
is a linear combination of the Weierstrass 
$\wp$-function and its shift by three three 
half periods. This is a non-autonomous system, 
because the Hamiltonian depends on the ``time'' 
$\tau$ through the $\tau$-dependence of the 
$\wp$-function. 

Levin and Olshanetsky \cite{bib:Le-Ol} 
pointed out that Manin's equation resembles 
the so called Calogero-Moser systems, i.e., 
the various extensions \cite{bib:Ol-Pe-review} 
of the integrable many-body systems first 
discovered by Calogero \cite{bib:Calogero}.  
More precisely, the Hamiltonian $\calH$ is 
identical to a special case (the rank-one 
elliptic model) of Inozemtsev's extensions 
\cite{bib:In-Me,bib:Inozemtsev} of the 
Calogero-Moser systems.  Levin and Olshanetsky 
called this relation the ``Painlev\'e-Calogero 
correspondence''.  

One will naturally ask if this correspondence 
can be extended to the other Painlev\'e 
equations.  Manin himself raised this problem 
in his paper.  Olshanetsky \cite{bib:Ol-review} 
conjectured that a degenerate version of 
Inozemtsev's elliptic model will emerge therein.  

This paper aims to answer this question affirmatively. 
A guiding principle is the degeneration relation 
of the six Painlev\'e equations \cite{bib:Ok-garnier}. 
This relation can be schematically expressed as follows: 
\beq
  \begin{array}{ccccccc}
    \PVI & \longrightarrow & \PV & \longrightarrow & 
    \PIV & & \\
         &                 & \downarrow & & 
    \downarrow & & \\
         &                 & \PIII & \longrightarrow & 
    \PII & \longrightarrow & \PI 
  \end{array}
  \nonumber 
\eeq
This diagram means, for instance, that $\PV$ can be 
derived from $\PVI$ by a degeneration process, which 
amounts to confluence of singular points of the 
aforementioned linear ordinary differential equation 
in the isomonodromic approach.  We shall trace this 
process carefully on the ``Calogero side'', and find 
a $\PV$-version of Manin's equation.  In principle, 
one can thus find an analogue of Manin's equation for 
all the six Painlev\'e equations (though, actually, 
one can resort to a more direct approach that bypasses 
the complicated degeneration process).  

Remarkably (or rather naturally?), all the six equations 
on the Calogero side turn out to become a (non-autonomous) 
Hamiltonian system with a Hamiltonian of the normal form 
$\calH = p^2/2 + V(q)$.  Furthermore, the Hamiltonians 
on the Calogero side of $\PV$ and $\PIV$ coincide with 
the Hamiltonians of the (rank one) hyperbolic and 
rational models in Inozemtsev's classification 
\cite{bib:In-Me} (which were discovered by Levi and 
Wojciechowski \cite{bib:Le-Wo} before Inozemtsev's 
work).  Those corresponding to the other three Painlev\'e 
equations are not included therein, but may be thought of 
as a further degeneration of the hyperbolic and rational 
models.  

One can further proceed to the higher rank models, 
and ask if there is still a Painlev\'e-Calogero 
correspondence.  We shall show that this is 
also the case.  The Painlev\'e side of the 
correspondence is a kind of multi-dimensional 
extensions of the Painlev\'e equations.  
They are obviously different from another 
multi-dimensional extension called the 
``Garnier systems'' \cite{bib:Ok-garnier}.  
For this reason, we call our multi-dimensional 
extension a {\it multi-component} version of 
the Painlev\'e equations.  

This paper is organized as follows.  Section 2 is 
a brief review of the work of Fuchs, Painlev\'e 
and Manin.  Section 3 deals with $\PV$, $\PIV$ and 
$\PIII$.  The degeneration process is discussed 
in detail for the case of $\PV$.  The direct 
approach is illustrated for the case of $\PIV$ and 
$\PIII$. Section 4 shows a reformulation of the 
foregoing calculations in the Hamiltonian formalism.  
The status of $\PII$ and $\PI$ is also clarified 
therein.  Section 5 is devoted to the higher rank 
Inozemtsev Hamiltonians and the multi-component 
Painlev\'e equations.  Section 6 is for concluding 
remarks.  Part of technical details are gathered 
in Appendices.

\section{Painlev\'e-Calogero Correspondence for $\PVI$}

We here briefly review the work of Fuchs, Painlev\'e 
and Manin. 

Fuchs rewrites $\PVI$ into the following form: 
\beq
    && t(1 - t)\calL_t \int_\infty^\lambda
       \frac{dz}{\sqrt{z(z - 1)(z - t)}}
    \nonumber \\
    &=& \sqrt{\lambda(\lambda - 1)(\lambda - t)} 
        \left[ \alpha + \frac{\beta t}{\lambda^2} 
        + \frac{\gamma(t - 1)}{(\lambda - 1)^2} 
        + \left(\delta - \frac{1}{2}\right) 
          \frac{t(t-1)}{(\lambda - t)^2} \right]. 
    \label{eq:PVI-Fuchs}
\eeq
Here $\calL_t$ is the linear differential operator 
(Picard-Fuchs operator) 
\beq
    \calL_t = t(1 - t) \frac{d^2}{dt^2} 
    + (1 - 2t) \frac{d}{dt} - \frac{1}{4}, 
\eeq
which also appears in the Picard-Fuchs equation 
of  complete elliptic integrals.  In this respect, 
$\PVI$ may be thought of as an inhomogeneous 
(and nonlinear) analogue of the Picard-Fuchs 
equation.  

Painlev\'e and Manin make use of a parametrization 
of the elliptic curve 
\beq
    y^2 = z(z - 1)(z - t) 
\eeq
by the Weierstrass $\wp$-function.  Let $\wp(u)$ 
be the $\wp$-function with primitive periods 
$1$ and $\tau$: 

\beq
    \wp(u) = \wp(u \mid 1,\tau) 
    = \frac{1}{u^2} 
      + \sum_{(m,n) \not= (0,0)} 
        \left( \frac{1}{(u + m + n\tau)^2} 
        - \frac{1}{(m + n\tau)^2} \right), 
\eeq
The parametrization is now given by 
\beq
    z = \frac{\wp(u) - e_1}{e_2 - e_1}, \quad 
    y = \frac{\wp'(u)}{2(e_2 - e_1)^{3/2}}, 
\eeq
where $e_n = \wp(\omega_n)$, $n = 1,2,3$ are 
the values of $\wp(u)$ at the three half period 
points $\omega_1 = 1/2$, $\omega_2 = - (1+\tau)/2$, 
$\omega_3 = \tau/2$. 

Manin's excellent idea is to do a simultaneous 
change of the dependent variable $\lambda \to q$ by 
\beq
    \lambda = \frac{\wp(q) - e_1}{e_2 - e_1}, 
    \label{eq:PVI-lam-q}
\eeq
and the independent vrariable $t \to \tau$ by 
\beq
    t = \frac{e_3 - e_1}{e_2 - e_1}. 
    \label{eq:PVI-t-tau}
\eeq
Manin presents the beautiful formula 
\beq
    \frac{d\tau}{dt} 
    = \frac{\pi i}{t(t - 1)(e_2 - e_1)}, 
    \label{eq:PVI-dtau/dt}
\eeq
for the Jacobian of the latter, which plays 
a key role in his calculations.  $\PVI$ is 
thereby mapped to the equation 
\beq
    (2 \pi i)^2 \frac{d^2 q}{d\tau^2} 
    = \sum_{n=0}^3 \alpha_n \wp'(q + \omega_n), 
    \label{eq:PVI-Manin}
\eeq
where the parameters on the right hand side 
are connected with the parameters of $\PVI$ 
as $\alpha_0 = \alpha$, $\alpha_1 = -\beta$, 
$\alpha_2 = \gamma$, $\alpha_3 = -\delta + 1/2$. 
This equation is equivalent to the Hamiltonian 
system 
\beq
    2\pi i\frac{dq}{d\tau} = \frac{\rd\calH}{\rd p}, 
    \quad 
    2\pi i\frac{dp}{d\tau} = - \frac{\rd\calH}{\rd q} 
\eeq
with the Hamiltonian 
\beq
    \calH = \frac{p^2}{2} 
    - \sum_{n=0}^3 \alpha_n \wp(q + \omega_n). 
\eeq

\section{Correspondence for $\PV$, $\PIV$ and $\PIII$} 

\subsection{Degeneration of $\PVI$ to $\PV$}

The degeneration of $\PVI$ to $\PV$ is achieved by 
rescaling the time variable and the parameters as 
\beq
    t = 1 + \epsilon \tilde{t}, \quad 
    \alpha = \tilde{\alpha}, \quad 
    \beta = \tilde{\beta}, \quad 
    \gamma = \frac{\tilde{\gamma}}{\epsilon} 
      - \frac{\tilde{\delta}}{\epsilon^2},  \quad 
    \delta = \frac{\tilde{\delta}}{\epsilon^2} 
\eeq
and letting $\epsilon \to 0$ while leaving 
$\tilde{\alpha},\ldots,\tilde{\gamma}$ and $\tilde{t}$ 
finite \cite{bib:Ok-garnier}.  

The building blocks of Fuchs' equation (\ref{eq:PVI-Fuchs}) 
turn out to survive this scaling limit as follows: 
\begin{enumerate}
\item The Picard-Fuchs operator: 
\beq
    t(1 - t)\calL_t \ \longrightarrow \ 
    \tilde{t}^2 \frac{d^2}{d\tilde{t}^2} 
    + \tilde{t} \frac{d}{d\tilde{t}} 
    = \left(\tilde{t}\frac{d}{d\tilde{t}}\right)^2. 
    \nonumber 
\eeq
\item The sum $\alpha + \cdots$ of four terms 
on the right hand side: 
\beq
   \alpha + \frac{\beta t}{\lambda^2} 
    + \frac{\gamma(t - 1)}{(\lambda - 1)^2} 
    + \left(\delta - \frac{1}{2}\right)
      \frac{t(t - 1)}{(\lambda - t)^2} 
   \ \longrightarrow \ 
   \tilde{\alpha} + \frac{\tilde{\beta}}{\lambda^2} 
    + \frac{\tilde{\gamma} \tilde{t}}{(\lambda - 1)^2} 
    + \frac{\tilde{\delta} \tilde{t}^2 (\lambda + 1)}
           {(\lambda - 1)^3}. 
   \nonumber 
\eeq
\item The quare root on the right hand side: 
\beq
    \sqrt{\lambda(\lambda - 1)(\lambda - t)} 
    \ \longrightarrow \ 
    \sqrt{\lambda}(\lambda - 1). 
    \nonumber 
\eeq
\item The incomplete elliptic integral: 
\beq
    \int_\infty^\lambda \frac{dz}{\sqrt{z(z-1)(z-t)}} 
    \ \longrightarrow \ 
    \int_\infty^\lambda \frac{dz}{\sqrt{z}(z-1)}. 
    \nonumber 
\eeq
\end{enumerate}
In particular, the degeneration of $\PVI$ 
to $\PV$ is associated with the degeneration of 
the elliptic curve to a rational curve, 
\beq
    y^2 = z(z - 1)(z - t) \quad 
    \longrightarrow \quad 
    y^2 = z(z - 1)^2, 
\eeq
or, equivalently, the degeneration of the torus 
$\bbC/(\bbZ + \tau\bbZ)$ to the cylinder 
$\bbC/\bbZ$. 

Thus, rewriting $\tilde{\alpha},\tilde{\beta},
\tilde{\gamma},\tilde{\delta}$ and $\tilde{t}$ 
to $\alpha,\beta,\gamma,\delta$ and $t$, 
we obtain the following equation as a $\PV$-version 
of Fuchs' equation: 
\beq
    \left(t\frac{d}{dt}\right)^2 
    \int_\infty^\lambda \frac{dz}{\sqrt{z}(z-1)} 
    = \sqrt{\lambda}(\lambda - 1) 
      \left(\alpha + \frac{\beta}{\lambda^2} 
      + \frac{\gamma t}{(\lambda-1)^2} 
      + \frac{\delta t^2(\lambda+1)}{(\lambda-1)^3}\right). 
    \label{eq:PV-Fuchs}
\eeq

\subsection{Analogue of Manin's equation for $\PV$} 

As an counterpart of the $q$-variable for $\PVI$, 
we now consider 
\beq
    q = \int_\infty^\lambda \frac{dz}{\sqrt{z}(z - 1)}. 
\eeq
If one prefers to being more faithful to Manin's 
parametrization, one should rather define $q$ as 
\beq
    q = \frac{1}{2\pi i}\int_\infty^\lambda 
        \frac{dz}{\sqrt{z}(z - 1)}, 
    \nonumber 
\eeq
because $2(e_2 - e_1)^{1/2} \to 2\pi i$ as 
$\Im\tau \to +\infty$ (see Appendix B).  
Since there is no substantial difference, 
let us take the first definition that is 
slightly simpler for calculations.  

Let us rewrite (\ref{eq:PV-Fuchs}) in terms of  $q$. 
The integral can be readily calculated as 
\beq
    q = \log\left(\frac{\sqrt{\lambda}-1}
        {\sqrt{\lambda}+1} \right), 
\eeq
so that the inverse relation can be written 
\beq
    \sqrt{\lambda} = - \coth(q/2). 
    \label{eq:PV-lam-q}
\eeq
Terms on the right hand side of (\ref{eq:PV-Fuchs}) 
can be calculated as follows: 
\beq
    \sqrt{\lambda} (\lambda - 1)
    &=& - \frac{\cosh(q/2)}{\sinh^3(q/2)}, 
    \nonumber \\
    \sqrt{\lambda} (\lambda - 1)\frac{1}{\lambda^2} 
    &=& - \frac{\sinh(q/2)}{\cosh^3(q/2)}, 
    \nonumber \\
    \sqrt{\lambda} (\lambda - 1)\frac{1}{(\lambda - 1)^2} 
    &=& - \frac{1}{2} \sinh(q), 
    \nonumber \\ 
    \sqrt{\lambda} (\lambda - 1)
    \frac{(\lambda + 1)}{(\lambda - 1)^3} 
    &=& - \frac{\lambda^{3/2} + \lambda^{1/2}}
               {(\lambda - 1)^2} 
      = - \frac{1}{4}\sinh(2q). 
    \nonumber 
\eeq
The differential equation for $q$ eventually takes 
the form 
\beq
    \left(t\frac{d}{dt}\right)^2 q 
    = - \frac{\rd V(q)}{\rd q}, 
    \label{eq:PV-Manin}
\eeq
where 
\beq
    V(q) = - \frac{\alpha}{\sinh^2(q/2)} 
           - \frac{\beta}{\cosh^2(q/2)} 
           + \frac{\gamma t}{2} \cosh(q) 
           + \frac{\delta t^2}{8} \cosh(2q). 
\eeq
This gives a $\PV$-version of Manin's equation.  
Note that this equation can be readily converted 
to a Hamiltonian system with the Hamiltonian 
$\calH = p^2/2 + V(q)$.  

\remark

A very similar change of dependent variable for $\PV$ is 
discussed in the book of Iwasaki et al. \cite{bib:IKSY}.

\subsection{Idea of direct approach}

Although the degeneration process can be 
continued to the other Painlev\'e equations, 
we now present a more direct approach.  
Note that the integrand is connected with 
the coefficient of $(d\lambda/dt)^2$ in 
the original Painlev\'e equation by the 
following very simple relation: 
\beq
    \frac{1}{\sqrt{z(z - 1)(z - t)}} 
    &=& \exp\left[ -\int \frac{1}{2}
        \left(\frac{1}{z} + \frac{1}{z-1} 
        + \frac{1}{z-t}\right)dz \right], 
    \nonumber \\
    \frac{1}{\sqrt{z}(z - 1)} 
    &=& \exp\left[ -\int 
        \left(\frac{1}{2z} + \frac{1}{z-1} 
        \right)dz \right]. 
    \nonumber 
\eeq
If this is a correct prescription, one 
will be able to define the $q$-variable 
for $\PIII$ and $\PII$ directly without 
the cumbersome degeneration process.  
This is indeed the case, as we shall show 
below.

\subsection{$q$-variable for $\PIV$} 

Since the expected integrand is given by 
\beq
    \exp\left( -\int \frac{dz}{2z} \right) 
    = \frac{1}{\sqrt{z}}, 
\eeq
we define 
\beq
    q = \int^\lambda \frac{dz}{\sqrt{z}} 
      = 2\sqrt{\lambda}. 
\eeq
This can be solved for $\lambda$ as 
\beq
    \lambda = \left(\frac{q}{2}\right)^2.  
    \label{eq:PIV-lam-q} 
\eeq
Honest calculations show that all derivative 
terms of $\PIV$ can be absorbed by the 
second derivative of $q$: 
\beq
    \frac{d^2q}{dt^2} 
    &=& \frac{1}{\sqrt{\lambda}} \frac{d^2\lambda}{dt^2} 
      - \frac{1}{2\lambda\sqrt{\lambda}} \left( 
        \frac{d\lambda}{dt} \right)^2 
    \nonumber \\
    &=& \frac{1}{\sqrt{\lambda}} \left( 
        \frac{3}{2} \lambda^3 + 4t\lambda^2 
        + 2(t^2 - \alpha)\lambda 
        + \frac{\beta}{\lambda} \right). 
\eeq
Substituting $\lambda = (q/2)^2$ gives the second 
order differential equation 
\beq
    \frac{d^2q}{dt^2} = - \frac{\rd V(q)}{\rd q} 
    \label{eq:PIV-Manin} 
\eeq
with the potential 
\beq
    V(q) = - \frac{1}{2}\left(\frac{q}{2}\right)^6 
           - 2t \left(\frac{q}{2}\right)^4 
           - 2(t^2 - \alpha)\left(\frac{q}{2}\right)^2 
           + \beta \left(\frac{q}{2}\right)^{-2}. 
\eeq

\subsection{$q$-variable for $\PIII$} 

The integrand is expected to be given by 
\beq
    \exp\left( - \int \frac{dz}{z} \right) 
    = \frac{1}{z}. 
\eeq
We consider 
\beq
    q = \int^\lambda \frac{dz}{z} = \log \lambda 
\eeq
and its inversion 
\beq
    \lambda = e^{q}. 
    \label{eq:PIII-lam-q} 
\eeq
All derivatives terms of $\PIII$  are now absorbed 
by the second derivative of $q$ with respect to 
$\log t$: 
\beq
    \left(t\frac{d}{dt}\right)^2 q 
    &=& \frac{t^2}{\lambda}\frac{d^2\lambda}{dt^2} 
      + \frac{t}{\lambda}\frac{d\lambda}{dt} 
      - \frac{t^2}{\lambda^2}
        \left(\frac{d\lambda}{dt}\right)^2 
    \nonumber \\
    &=& \frac{\alpha\lambda}{4} 
      + \frac{\beta t}{4\lambda} 
      + \frac{\gamma\lambda^2}{4} 
      + \frac{\delta t^2}{4\lambda^2}. 
\eeq
Substituting $\lambda = e^q$ gives the second 
order equation 
\beq
    \left(t\frac{d}{dt}\right)^2q 
    = - \frac{\rd V(q)}{\rd q} 
    \label{eq:PIII-Manin} 
\eeq
with the potential 
\beq
    V(q) = - \frac{\alpha}{4} e^q 
           + \frac{\beta t}{4} e^{-q} 
           - \frac{\gamma}{8} e^{2q} 
           + \frac{\delta t^2}{8} e^{-2q}. 
\eeq

\subsection{Summary}

Let us summarize the results of this section.  

\begin{theorem} 
The foregoing change of variable $\lambda \to q$ 
maps $\PV$, $\PIV$ and $\PIII$ to a second order 
differential equation for the new dependent 
variable $q$. These equations are equivalent 
to a non-autonomous Hamiltonian system with 
a Hamiltonian of the normal form 
$\calH = p^2/2 + V(q)$: 
\newline
{\rm ($\PV$)} 
The Hamiltonian system takes the form 
\beq
    t\frac{dq}{dt} = \frac{\rd\calH}{\rd p}, \quad 
    t\frac{dp}{dt} = - \frac{\rd\calH}{\rd q} 
\eeq
with the Hamiltonian 
\beq
    \calH = \frac{p^2}{2} 
           - \frac{\alpha}{\sinh^2(q/2)} 
           - \frac{\beta}{\cosh^2(q/2)} 
           + \frac{\gamma t}{2} \cosh(q) 
           + \frac{\delta t^2}{8} \cosh(2q). 
\eeq
\newline
{\rm ($\PIV$)} 
The Hamiltonian system takes the form 
\beq
    \frac{dq}{dt} = \frac{\rd\calH}{\rd p}, \quad 
    \frac{dp}{dt} = - \frac{\rd\calH}{\rd q} 
\eeq
with the Hamiltonian 
\beq
    \calH = \frac{p^2}{2} 
           - \frac{1}{2}\left(\frac{q}{2}\right)^6 
           - 2t \left(\frac{q}{2}\right)^4 
           - 2(t^2 - \alpha)\left(\frac{q}{2}\right)^2 
           + \beta \left(\frac{q}{2}\right)^{-2}. 
\eeq
\newline
{\rm ($\PIII$)}
The Hamiltonian system takes the form 
\beq
   t\frac{dq}{dt} = \frac{\rd\calH}{\rd p}, \quad 
   t\frac{dp}{dt} = - \frac{\rd\calH}{\rd q} 
\eeq
with the Hamiltonian 
\beq
   \calH = \frac{p^2}{2} 
           - \frac{\alpha}{4} e^q 
           + \frac{\beta t}{4} e^{-q} 
           - \frac{\gamma}{8} e^{2q} 
           + \frac{\delta t^2}{8} e^{-2q}. 
\eeq
\end{theorem}

\remark
\begin{enumerate}
\item 
The Hamiltonians for $\PV$ and $\PIV$ coincide with 
those of the hyperbolic and rational models of 
Inozemtsev \cite{bib:In-Me}, Levi and Wojciechowski 
\cite{bib:Le-Wo}.  The Hamiltonian for $\PIII$ has 
no counterpart in their work, but nowadays can be 
found in the literature \cite{bib:vanDiejen}.  
\item  
The foregoing construction of the $q$-variable 
does not literally work for $\PII$ and $\PI$, 
because there is no $(d\lambda/dt)^2$ term.  
The status of these equations will be clarified 
in the next section from a different point of view. 
\end{enumerate}

\section{Hamiltonian formalism of correspondence}

\subsection{Hamiltonians of Painlev\'e equations} 

All the six Painlev\'e equations are known to be 
expressed in the Hamiltonian form 
\beq
    \frac{d\lambda}{dt} = \frac{\rd H}{\rd\mu}, \quad 
    \frac{d\lambda}{dt} = - \frac{\rd H}{\rd\lambda} 
    \nonumber 
\eeq
with a suitable choice of the canonical conjugate 
variable $\mu$ and the Hamiltonian $H$ \cite{bib:Malmquist}. 
This expression is by no means unique;  we here 
consider the following Hamiltonians \cite{bib:Ok-garnier}.  
These Hamiltonians are referred to as the ``polynomial 
Hamiltonians'' because they are  polynomials in $\lambda$ 
and $\mu$: 
\beq
    &(\PVI)& 
    H = \frac{\lambda(\lambda - 1)(\lambda - t)}{t(t - 1)} 
        \left[\mu^2 - \left(\frac{\kappa_0}{\lambda} 
          + \frac{\kappa_1}{\lambda - 1} 
          + \frac{\theta - 1}{\lambda - t}
            \right) \mu 
          + \frac{\kappa}{\lambda(\lambda - 1)}\right]. 
    \nonumber \\
    &(\PV)& 
    H = \frac{\lambda(\lambda - 1)^2}{t} 
        \left[\mu^2 - \left(\frac{\kappa_0}{\lambda} 
          + \frac{\theta_1}{\lambda - 1} 
          - \frac{\eta_1 t}{(\lambda - 1)^2} 
          \right) \mu 
          + \frac{\kappa}{\lambda(\lambda - 1)} \right]. 
    \nonumber \\
    &(\PIV)& 
    H = 2\lambda \left[\mu^2 - \left(\frac{\lambda}{2} 
          + t + \frac{\kappa_0}{\lambda}\right) \mu 
          + \frac{\theta_\infty}{2} \right]. 
    \nonumber \\
    &(\PIII)&
    H = \frac{\lambda^2}{t}
        \left[\mu^2 - \left(\eta_\infty 
          + \frac{\theta_0}{\lambda} 
          - \frac{\eta_0 t}{\lambda^2}\right) \mu 
          + \frac{\eta_\infty(\theta_0 + \theta_\infty)} 
                 {2\lambda} 
          \right]. 
    \nonumber \\
    &(\PII)& 
    H = \frac{\mu^2}{2} 
        - \left(\lambda^2 + \frac{t}{2}\right) \mu 
        - \left(\alpha + \frac{1}{2}\right) \lambda. 
    \nonumber \\
    &(\PI)& 
    H = \frac{\mu^2}{2} - 2\lambda^3 - t\lambda. 
    \nonumber 
\eeq
Here $\kappa_0,\kappa_1,\theta$, etc. are constants 
that are connected with the parameters $\alpha,\beta,
\gamma,\delta$ of the Painlev\'e equations by simple 
algebraic relations: 
\beq
    &(\PVI)& 
    \alpha = \frac{(\kappa_0 + \kappa_1 + \theta - 1)^2}{2} 
           - 2\kappa , \ 
    \beta = - \frac{\kappa_0^2}{2}, \ 
    \gamma = \frac{\kappa_1^2}{2}, \ 
    \delta = \frac{1 - \theta^2}{2}, \ 
    \nonumber \\
    &(\PV)& 
    \alpha = \frac{(\kappa_0 + \theta_1)^2}{2} 
           - 2\kappa, \ 
    \beta = - \frac{\kappa_0^2}{2}, \ 
    \gamma = \eta_1(\theta_1 + 1), \ 
    \delta = - \frac{\eta_1^2}{2}. 
    \nonumber \\
    &(\PIV)& 
    \alpha = 2\theta_\infty - \kappa_0 + 1, \ 
    \beta = - 2\kappa_0^2. 
    \nonumber \\
    &(\PIII)& 
    \alpha = - 4\eta_\infty\theta_\infty, \ 
    \beta = 4\eta_0(\theta_0 + 1), \ 
    \gamma = 4\eta_\infty^2, \ 
    \delta = - 4\eta_0^2. 
    \nonumber 
\eeq

\subsection{How to find canonical transformations} 

The goal of this section is to show that the 
Painlev\'e-Calogero correspondence is, in fact, 
a (time-dependent) canonical transformation of 
two Hamiltonian systems.  By this, we mean that 
the functional relation between $\lambda$ and $q$ 
can be extended to $(\lambda,\mu)$ and $(q,p)$ 
so as to satisfy the equation 
\beq
    \mu d\lambda - H dt 
    = \const \cdot (pdq - \calH dT) 
    + \mbox{\rm exact form}. 
\eeq
with a suitably redefined time variable $T$ 
(such as the logarithmic time $\log t$ in $\PV$ 
and $\PIII$). The constant factor on the right 
hand side is inserted simply for convenience; 
if necessary, one can normalize the constant 
to $1$ by suitably rescaling $p,q,\calH$ and 
$T$.  For this reason, wel call this type of 
coordinate transformation a ``canonical'' 
transformation even if the constant factor is 
not equal to $1$.  

Let us illustrate, in the case of $\PVI$, how to 
find such a canonical transformation.  Suppose 
that $\lambda$ and $\mu$ be a solution of $\PVI$ in 
the aforementioned Hamiltonian formalism, and that 
$q$ be a corresponding solution of Manin's equation.  
The canonical equation for $\lambda$ takes the form 
\beq
    \frac{d\lambda}{dt} 
    = \frac{\lambda(\lambda - 1)(\lambda - t)}{t(t - 1)} 
      \left(2\mu - \frac{\kappa_0}{\lambda} 
        - \frac{\kappa_1}{\lambda - 1} 
        - \frac{\theta - 1}{\lambda - t} \right). 
    \nonumber 
\eeq
This equations can be solve for $\mu$: 
\beq
    \mu = \frac{t(t - 1)}{2\lambda(\lambda - 1)(\lambda - t)} 
          \frac{d\lambda}{dt} 
        + \frac{1}{2} \left(\frac{\kappa_0}{\lambda} 
          + \frac{\kappa_1}{\lambda - 1}  
          + \frac{\theta - 1}{\lambda - t} \right). 
    \nonumber 
\eeq
Our task is to rewrite the right hand side in terms 
of $p$ and $q$.  We first consider $d\lambda/dt$. 
Differentiating (\ref{eq:PVI-lam-q}) against $t$ gives 
\beq
    \frac{d\lambda}{dt} 
    = \left(\frac{\wp'(q)}{e_2 - e_1}\frac{dq}{d\tau} 
      + f_\tau(q)\right) \frac{d\tau}{dt}, 
    \nonumber 
\eeq
where we have introduced the functions 
\beq
    f(u) = \frac{\wp(u) - e_1}{e_2 - e_1}, \quad 
    f_\tau(u) = \frac{\rd f(u)}{\rd \tau}. 
\eeq
The derivative $dq/d\tau$ can be read off from 
the canonical equation for $q$: 
\beq
    \frac{dq}{d\tau} 
    = \frac{1}{2\pi i}\frac{\rd\calH}{\rd p} 
    = \frac{p}{2\pi i}. 
    \nonumber 
\eeq
As for the Jacobian $d\tau/dt$, Manin's formula 
(\ref{eq:PVI-dtau/dt}) is available.  One can 
thus express $d\lambda/dt$ as a function of 
$p,q$ and $\tau$.  The other part of the 
foregoing expression of $\mu$ contains $\lambda$ 
only, which can be readily converted to a 
function of $q$ and $\tau$ by (\ref{eq:PVI-lam-q}). 
We thus obtain the following expression of $\mu$: 
\beq
    \mu &=& \frac{e_2 - e_1}{\wp'(q)} p 
    + \frac{2\pi i(e_2 - e_1)^2}{\wp'(q)^2} f_\tau(q) 
    \nonumber \\
    && + \frac{e_2 -e_1}{2} \left( 
         \frac{\kappa_0}{\wp(q) - e_1} 
         + \frac{\kappa_1}{\wp(q) - e_2} 
         + \frac{\theta - 1}{\wp(q) - e_3} \right). 
    \label{eq:PVI-mu-pq} 
\eeq         

We now move the point of view, and think of 
(\ref{eq:PVI-lam-q}) and (\ref{eq:PVI-mu-pq}) 
as defining a coordinate transformation 
$(\lambda,\mu) \to (q,p)$.  This gives a 
canonical transformation that we have 
sought for: 

\begin{theorem}\label{th:PVI-ct}
(\ref{eq:PVI-lam-q}) and (\ref{eq:PVI-mu-pq}) 
define a canonical transformation that connects 
the Hamiltonian form of $\PVI$ and Manin's 
Hamiltonian system. The canonical coordinates 
and the Hamiltonians of the two systems obey 
the equation 
\beq 
    \mu d\lambda - H dt 
    = pdq - \calH \frac{d\tau}{2\pi i}
    + \mbox{\rm exact form}. 
\eeq
\end{theorem}

\subsection{Proof of Theorem \ref{th:PVI-ct}}

Total differential of  (\ref{eq:PVI-lam-q}) 
gives 
\beq
    d\lambda = \frac{\wp'(q)}{e_2 - e_1}dq 
    + f_\tau(q)d\tau, 
    \nonumber 
\eeq
so that $\mu d\lambda$ can be expressed as 
\beq
   \mu d\lambda 
   &=& \left( \frac{e_2 - e_1}{\wp'(q)} p 
       + \frac{2\pi i(e_2- e_1)^2}{\wp'(q)^2} 
         f_\tau(q) \right) 
       \left( \frac{\wp'(q)}{e_2 - e_1} dq 
       + f_\tau(q) d\tau \right) 
   \nonumber \\
   && + \frac{1}{2}\left(\frac{\kappa_0}{\lambda} 
        + \frac{\kappa_1}{\lambda - 1} 
        + \frac{\theta - 1}{\lambda - t} \right) 
        d\lambda 
   \nonumber \\
   &=& pdq + (\mathrm{A}) + (\mathrm{B}) 
       + (\mathrm{C}), 
   \nonumber 
\eeq
where 
\beq
    (\mathrm{A}) 
    &=& \frac{2\pi i(e_2 -e_1)}{\wp'(q)} 
        f_\tau(q) dq, 
    \nonumber \\
    (\mathrm{B}) 
    &=& \left( \frac{e_2 - e_1}{\wp'(q)} p 
        + \frac{2\pi i(e_2 - e_1)^2}{\wp'(q)^2} 
          f_\tau(q) \right)
        f_\tau(q) d\tau, 
    \nonumber \\
    (\mathrm{C}) 
    &=& \frac{1}{2}\left( \frac{\kappa_0}{\lambda} 
        + \frac{\kappa_1}{\lambda - 1} 
        + \frac{\theta - 1}{\lambda - t} \right) 
        d\lambda.
    \nonumber 
\eeq
As we shall prove in Appendix A, (A) can be 
further rewritten 
\beq
    (\mathrm{A}) 
    = \left[ \frac{\wp(q + \omega_3)}{4\pi i} 
      - \pi\left(\frac{f_\tau(q)}{f'(q)}\right)^2 
      \right] d\tau 
      + \mbox{\rm exact form}, 
    \label{eq:term-A}
\eeq
where $f'(u)$ denotes the $u$-derivative 
of $f(u)$: 
\beq
    f'(u) = \frac{\rd f(u)}{\rd u} 
    = \frac{\wp'(u)}{e_2 - e_1}. 
\eeq
For (B) and (C), we have 
\beq
    (\mathrm{B}) 
    &=& \left[ \frac{f_\tau(q)}{f'(q)} p 
        + 2\pi i\left(\frac{f_\tau(q)}{f'(q)}
        \right)^2 \right] d\tau, 
    \nonumber \\
    (\mathrm{C}) 
    &=& \frac{\theta - 1}{2(\lambda - t)} dt 
      + \frac{1}{2}\left( 
        \kappa_0\log\lambda 
        + \kappa_1\log(\lambda - 1) 
        + (\theta - 1)\log(\lambda - t) \right) 
    \nonumber \\
    &=& \frac{\theta - 1}{2(\lambda - t)} dt 
      + \mbox{\rm exact form}. 
    \nonumber 
\eeq
Thus we find that 
\beq
    \mu d\lambda - Hdt 
    = pdq - \tilde{\calH} \frac{d\tau}{2\pi i} 
    + \mbox{\rm exact form}, 
\eeq
where 
\beq
    \tilde{\calH} 
    = 2\pi i\frac{dt}{d\tau} \left(
        H - \frac{\theta - 1}{2(\lambda - t)} \right) 
    - 2\pi i\left[ 
        \frac{\wp(q + \omega_3)}{4\pi i} 
        + \frac{f_\tau(q)}{f'(q)} p 
        + \pi i\left(\frac{f_\tau(q)}{f'(q)}\right)^2 
      \right]. 
\eeq
Our task is to prove that the transformed Hamiltonian 
$\tilde{\calH}$ coincides, modulo irrelevant terms, 
with the Hamiltonian of Manin's equation.  Here 
``irrelevant'' means that the term is a function 
of $t$ only.  Such a ``non-dynamical'' term can be 
absorbed by the ``exact form'' part of the foregoing 
relation of $1$-forms, thereby being negligible. 

Let us evaluate the contribution of $2\pi i(dt/d\tau) H$. 
By Manin's formula (\ref{eq:PVI-dtau/dt}) of 
$d\tau/dt$, and also by the identity 
\beq
    \lambda(\lambda - 1)(\lambda - t) 
    = \frac{\wp'(q)^2}{4(e_2 - e_1)^3}, 
    \nonumber
\eeq
we can rewrite $2\pi i(dt/d\tau)H$ as follows: 
\beq
    2\pi i\frac{dt}{d\tau}H 
    &=& \frac{\wp'(q)^2}{2(e_2 - e_1)^2}
        \left[ 
          \mu^2 
          - \left(\frac{\kappa_0}{\lambda} 
            + \frac{\kappa_1}{\lambda - 1} 
            + \frac{\theta - 1}{\lambda - t}
            \right) \mu 
          + \frac{\kappa}{\lambda(\lambda - 1)} 
        \right]
    \nonumber \\
    &=& \frac{\wp'(q)^2}{2(e_2 - e_1)^2} 
        \left[ 
          \mu 
          - \frac{1}{2}\left(
            \frac{\kappa_0}{\lambda} 
            + \frac{\kappa_1}{\lambda -1} 
            + \frac{\theta - 1}{\lambda - t} 
            \right)
        \right]^2 
    \nonumber \\
    && + \frac{\wp'(q)^2}{2(e_2 - e_1)^2} 
         \left[ 
           - \frac{1}{4}\left( 
             \frac{\kappa_0}{\lambda} 
             + \frac{\kappa_1}{\lambda -1} 
             + \frac{\theta - 1}{\lambda - t}
             \right)^2 
           + \frac{\kappa}{\lambda(\lambda - 1)} 
         \right]. 
    \nonumber 
\eeq
The first term on the right hand side is equal to 
\beq
    \frac{1}{2}\left(p + 2\pi i\frac{f_\tau(q)}{f'(q)}
    \right)^2 
    = \frac{p^2}{2} + 2\pi i\frac{f_\tau(q)}{f'(q)} p 
    + \left(2\pi i\frac{f_\tau(q)}{f'(q)}\right)^2, 
    \nonumber 
\eeq
by which the terms proportional to $f_\tau(q)/f'(q)$ 
and its square in the definition of $\tilde{\calH}$ are 
cancelled out.  The transformed Hamiltonian $\tilde{\calH}$ 
can now be expressed as 
\beq
    \tilde{\calH} 
    &=& \frac{p^2}{2} 
      - \frac{\wp'(q)^2}{2(e_2 - e_1)^2} 
      - \frac{(\theta - 1)t(t -1)(e_2 - e_1)}{\lambda -t}
    \nonumber \\
    && + \frac{\wp'(q)}{2(e_2 - e_1)^2} 
         \left[ 
           - \frac{1}{4}\left( 
             \frac{\kappa_0}{\lambda} 
             + \frac{\kappa_1}{\lambda -1} 
             + \frac{\theta - 1}{\lambda - t}
             \right)^2 
           + \frac{\kappa}{\lambda(\lambda - 1)} 
         \right]. 
\eeq
Note that this is already of the normal form 
$p^2/2 + \tilde{V}(q)$ with the potential 
\beq
    \tilde{V}(q) &=& 
       - \frac{\wp'(q)^2}{2(e_2 - e_1)^2} 
       - \frac{(\theta - 1)t(t -1)(e_2 - e_1)}{\lambda -t}
    \nonumber \\
    && + \frac{\wp'(q)}{2(e_2 - e_1)^2} 
         \left[ 
           - \frac{1}{4}\left( 
             \frac{\kappa_0}{\lambda} 
             + \frac{\kappa_1}{\lambda -1} 
             + \frac{\theta - 1}{\lambda - t}
             \right)^2 
           + \frac{\kappa}{\lambda(\lambda - 1)} 
         \right].  
\eeq

What remains is to express $\tilde{V}(q)$ as an
explicit function of $q$. To this end, we 
substitute the factor $\wp'(q)^2/2(e_2 - e_1)^2$ 
by $2(e_2 - e_1)\lambda(\lambda - 1)(\lambda - t)$, 
and rewrite the main part of $\tilde{V}(q)$ as a
linear combination of $\lambda$, $1/\lambda$, 
$1/(\lambda - 1)$ and $1/(\lambda - t)$.  
This leads to the following expression of 
$\tilde{V}(q)$: 
\beq
    \tilde{V}(q) 
    &=& 
    - \frac{(\kappa_0+\kappa_1+\theta-1)^2 - 4\kappa}{2} 
      (e_2 - e_1)\lambda 
    \nonumber \\
    && 
    - \frac{\kappa_0^2}{2} \cdot \frac{(e_2 - e_1)t}{\lambda} 
    - \frac{\kappa_1^2}{2} \cdot 
      \frac{(e_2 - e_1)(1 - t)}{\lambda - 1} 
    - \frac{(\theta - 1)^2 + 1}{2} \cdot 
      \frac{(e_2 - e_1)t(t - 1)}{\lambda - t} 
    \nonumber \\
    && 
    - \frac{1}{2} \wp(q + \omega_3) 
    + \mbox{\rm function of $t$ only}. 
    \nonumber 
\eeq
The final piece of the ring is the general formula 
\beq
    \wp(u + \omega_j) = e_j 
    + \frac{(e_j - e_k)(e_j - e_\ell)}{\wp(u) - e_j}
\eeq
where $(j,k,l)$ is a cyclic permutation of 
$(1,2,3)$.   This implies that 
\beq
    \frac{(e_2 - e_1)t}{\lambda} 
    &=& \wp(q + \omega_1) - e_1, 
    \nonumber \\
    \frac{(e_2 - e_1)(1 - t)}{\lambda - 1} 
    &=& \wp(q + \omega_2) - e_2, 
    \nonumber \\
    \frac{(e_2 - e_1)t(t - 1)}{\lambda - t} 
    &=& \wp(q + \omega_3) - e_3, 
    \nonumber 
\eeq
so that 
\beq
    \tilde{V}(q) 
    &=& 
    - \frac{(\kappa_0+\kappa_1+\theta-1)^2 - 4\kappa}{2} 
      \wp(q) 
    - \frac{\kappa_0^2}{2} \wp(q + \omega_1) 
    \nonumber \\
    && 
    - \frac{\kappa_1^2}{2} \wp(q + \omega_2) 
    - \frac{\theta^2}{2} \wp(q + \omega_3) 
    + \mbox{\rm function of $\tau$ only}. 
\eeq
Apart from the last term which is negligible, 
this potential is indeed the same as Manin's 
potential $V(q)$ (recall the algebraic relations 
connecting the constants $\kappa_0$, etc. 
and the parameters of $\PVI$).  This completes 
the proof of the theorem. 
\qed

\subsection{Canonical transformation for $\PV$}

This heuristic method for constructing a 
canonical transformation can be applied to 
the other Painlev\'e equations.  We here 
consider the case of $\PV$. 

Let $\lambda$ be a solution of $\PV$, 
$\mu$ the canonical conjugate variable, 
and $q$ the corresponding solution of 
(\ref{eq:PV-Manin}). The canonical equation 
for $\lambda$ can be written 
\beq
    \frac{d\lambda}{dt} 
    = \frac{\lambda(\lambda - 1)^2}{t} 
      \left(2\mu - \frac{\kappa_0}{\lambda} 
      - \frac{\theta_1}{\lambda - 1} 
      + \frac{\eta_1 t}{(\lambda - 1)^2}\right). 
    \nonumber 
\eeq
This equation can be solved for $\mu$ as 
\beq
    \mu = \frac{1}{2\lambda(\lambda - 1)^2} 
          t \frac{d\lambda}{dt} 
        + \frac{1}{2} \left(\frac{\kappa_0}{\lambda} 
          + \frac{\theta_1}{\lambda - 1} 
          - \frac{\eta_1 t}{(\lambda - 1)^2}\right). 
    \nonumber 
\eeq
By differentiating (\ref{eq:PV-lam-q}) against $t$ 
and using the canonical equation $tdq/dt = \rd\calH/\rd p 
= p$, we obtain the identity 
\beq
    t \frac{d\lambda}{dt} 
    = \sqrt{\lambda}(\lambda - 1) p, 
    \nonumber 
\eeq
which can be used to rewrite the expression of 
$\mu$ as 
\beq
    \mu &=& \frac{p}{2\sqrt{\lambda}(\lambda - 1)} 
        + \frac{1}{2}\left(\frac{\kappa_0}{\lambda} 
          + \frac{\theta_1}{\lambda - 1} 
          - \frac{\eta_1 t}{(\lambda - 1)^2}\right). 
    \label{eq:PV-mu-pq}
\eeq

We now reinterpret (\ref{eq:PV-lam-q}) and 
(\ref{eq:PV-mu-pq}) as defining a coordinate 
transformation $(\lambda,\mu) \to (q,p)$.  
This indeed turns out to give a canonical 
transformation that we have sought for: 

\begin{theorem} 
(\ref{eq:PV-lam-q}) and (\ref{eq:PV-mu-pq}) 
define a canonical transformation that connects 
$\PV$ and the $\PV$-version of Manin's Hamiltonian 
system.  The canonical coordinates and the 
Hamiltonians of the two systems obey the equation 
\beq
    \mu d\lambda - H dt 
    = \frac{1}{2}\left(pdq - \calH \frac{dt}{t}\right) 
    + \mbox{\rm exact form}. 
\eeq
\end{theorem}

\proof
Since $d\lambda$ and $dq$ are connected by 
the relation 
\beq
    d\lambda = \sqrt{\lambda}(\lambda - 1)dq, 
    \nonumber 
\eeq
$\mu d\lambda$ can be expressed as 
\beq
    \mu d\lambda 
    &=& \frac{1}{2} pdq 
      + \frac{1}{2}\left(
        \frac{\kappa_0}{\lambda} 
        + \frac{\theta_1}{\lambda - 1} 
        - \frac{\eta_1 t}{(\lambda - 1)^2}
        \right) d\lambda 
    \nonumber \\
    &=& \frac{1}{2}pdq 
      - \frac{\eta_1}{2(\lambda - 1)}dt 
      + \frac{1}{2}d\left(
        \kappa_0\log\lambda 
        + \theta_1\log(\lambda - 1) 
        + \frac{\eta_1 t}{\lambda - 1}\right), 
    \nonumber 
\eeq
so that 
\beq
    \mu d\lambda - Hdt 
    = \frac{1}{2}\left(pdq 
      - \tilde{\calH} \frac{dt}{t}\right) 
    + \mbox{\rm exact form}, 
\eeq
where
\beq
    \tilde{\calH} 
    = 2Ht + \frac{\eta_1 t}{\lambda - 1}. 
\eeq
We can rewrite $\tilde{\calH}$ to a normal 
form as 
\beq
    \tilde{\calH} 
    &=& 2\lambda(\lambda - 1)^2 \left[ 
          \mu - \frac{1}{2}\left( 
          \frac{\kappa_0}{\lambda} 
          + \frac{\theta_1}{\lambda - 1} 
          - \frac{\theta_1 t}{(\lambda - 1)^2}
          \right) 
        \right]^2 
    \nonumber \\
    && + 2\lambda(\lambda - 1)^2 \left[ 
          - \frac{1}{4}\left(
            \frac{\kappa_0}{\lambda} 
            + \frac{\theta_1}{\lambda - 1} 
            - \frac{\eta_1 t}{(\lambda -1)^2} 
            \right)^2 
          + \frac{\kappa}{\lambda(\lambda - 1)}
         \right] 
       + \frac{\eta_1 t}{\lambda -1} 
    \nonumber \\
    &=& \frac{p^2}{2} + \tilde{V}(q), 
\eeq
where 
\beq
    \tilde{V}(q) 
    &=& - \frac{\lambda(\lambda - 1)^2}{2} 
        \left(\frac{\kappa_0}{\lambda} 
          + \frac{\theta_1}{\lambda - 1} 
          - \frac{\eta_1 t}{(\lambda - 1)^2} 
        \right)^2 
      + 2\kappa(\lambda -1) 
      + \frac{\eta_1 t}{\lambda -1}. 
    \nonumber \\
    &=& - \left(\frac{\kappa_0}{2} + \frac{\theta_1^2}{2} 
           + \kappa_1 \theta_1  - 2\kappa\right) 
           \frac{1}{\sinh^2(q/2)} 
        + \frac{\kappa_0^2}{2} \frac{1}{\cosh^2(q/2)} 
    \nonumber \\
    &&  + \frac{\eta_1(\theta_1 + 1) t}{2} \cosh(q) 
        - \frac{\eta_1^2 t^2}{2} \cosh(2q) 
        + \mbox{\rm function of $t$ only}. 
\eeq
Apart from the last negligible term, 
this coincides with the potential $V(q)$ 
in the statement of the theorem. 
\qed

\subsection{Canonical transformation for $\PIV$} 

We now consider the case of $\PIV$. 

Let $\lambda$ be a solution of $\PIV$, 
$\mu$ the canonical conjugate variable, 
and $q$ the corresponding solution of 
(\ref{eq:PIV-Manin}).  The canonical equation 
for $\lambda$ can be written 
\beq
    \frac{d\lambda}{dt} 
    = 4 \lambda \mu 
    - (\lambda^2 + 2t \lambda  + 2\kappa_0), 
    \nonumber 
\eeq
which can be solved for $\mu$ as 
\beq
    \mu = \frac{1}{4\lambda}\frac{d\lambda}{dt} 
        + \frac{1}{4}\left(\lambda + 2t 
          + \frac{2 \kappa_0}{\lambda}\right). 
    \nonumber 
\eeq
By (\ref{eq:PIV-lam-q}) and the canonical equation 
$dq/dt = \rd\calH/\rd p = p$, we have the 
identity 
\beq
    \frac{d\lambda}{dt} 
    = \sqrt{\lambda} \frac{dq}{dt} 
    = \sqrt{\lambda} p, 
    \nonumber 
\eeq
so that 
\beq
    \mu = \frac{p}{4\sqrt{\lambda}} 
        + \frac{1}{4}\left(\lambda + 2t 
          + \frac{2 \kappa_0}{\lambda}\right). 
    \label{eq:PIV-mu-pq}
\eeq

\begin{theorem}
(\ref{eq:PIV-lam-q}) and (\ref{eq:PIV-mu-pq}) 
define a canonical transformation that connects 
$\PIV$ and the $\PIV$-version of Manin's 
Hamiltonian system.  The canonical coordinates 
and Hamiltonians of the two systems obey 
the equation 
\beq
    \mu d\lambda - H dt 
    = \frac{1}{4}(pdq - \calH dt) 
    + \mbox{\rm exact form}. 
\eeq
\end{theorem}

\proof
Since $d\lambda$ and $dq$ are connected by 
the relation 
\beq
    d\lambda = \sqrt{\lambda} dq, 
    \nonumber 
\eeq
$\mu d\lambda$ can be expressed as 
\beq
    \mu d\lambda 
    &=& \frac{1}{4} pdq 
      + \frac{1}{4}\left(\lambda + 2t 
          + \frac{2\kappa_0}{\lambda}\right)d\lambda 
    \nonumber \\
    &=& \frac{1}{4}pdq 
      - \frac{1}{2}\lambda dt 
      + \frac{1}{4}d\left(\frac{\lambda^2}{2} 
         + 2t\lambda + 2\kappa_0\log\lambda\right), 
    \nonumber 
\eeq
so that 
\beq
    \mu d\lambda - Hdt 
    = \frac{1}{4}(pdq - \tilde{\calH}dt) 
    + \mbox{\rm exact form}, 
\eeq
where 
\beq
    \tilde{\calH} = 4H + 2\lambda. 
\eeq
We can rewrite the transformed Hamiltonian 
$\tilde{\calH}$ to a normal form as 
\beq
    \tilde{\calH} 
    &=& 8\lambda\left[
           \mu - \frac{1}{2}\left( 
             \frac{\lambda}{2} + t 
             + \frac{\kappa_0}{\lambda}\right) 
        \right]^2 
      + 8\lambda\left[ 
           - \frac{1}{4}\left(
              \frac{\lambda}{2} + t 
              + \frac{\kappa_0}{\lambda} \right)^2 
           + \frac{\theta_\infty}{2} \right] 
       + 2\lambda 
    \nonumber \\
    &=& \frac{p^2}{2} + \tilde{V}(q), 
\eeq
where 
\beq
    \tilde{V}(q) 
    &=&  - 2\lambda\left(\frac{\lambda}{2} + t 
           + \frac{\kappa_0}{\lambda}\right)^2 
       + 4\theta_\infty\lambda 
       + 2\lambda
    \nonumber \\
    &=& - \frac{1}{2}\lambda^3 - 2t\lambda^2 
        - 2(t^2 + \kappa_0 - 2\theta_\infty - 1) \lambda 
        - 2\kappa_0^2 \lambda^{-1} 
    \nonumber \\
    && + \mbox{\rm function of $t$ only}. 
\eeq
Substituting $\lambda = (q/2)^2$ gives 
the potential $V(q)$ modulo an irrelevant term. 
\qed

\subsection{Canonical transformations for $\PIII$}

The situation of $\PIII$ is somewhat similar 
to $\PV$.  

Let $\lambda$, again, be a solution of $\PIII$, 
$\lambda$ the canonical conjuage variable, 
and $q$ be the corresponding solution of 
(\ref{eq:PIII-Manin}).  The canonical 
equation for $\lambda$ takes the form 
\beq
    \frac{d\lambda}{dt} 
    = \frac{\lambda^2}{t}
      \left(2\mu - \eta_\infty 
      - \frac{\theta_0}{\lambda} 
      + \frac{\eta_0 t}{\lambda^2} \right), 
    \nonumber 
\eeq 
which can be solved for $\mu$ as 
\beq
    \mu = \frac{t}{2\lambda^2} \frac{d\lambda}{dt} 
        + \frac{1}{2} \left(\eta_\infty 
          + \frac{\theta_0}{\lambda} 
          - \frac{\eta_0 t}{\lambda^2}\right). 
    \nonumber 
\eeq
By differentiating (\ref{eq:PIII-lam-q}) and 
using the canonical equation $tdq/dt = \rd\calH/\rd p 
= p$, the $t$-derivative of $\lambda$ can be 
written 
\beq
    t \frac{d\lambda}{dt} = \lambda p, 
    \nonumber 
\eeq
so that we obtain 
\beq
    \mu = \frac{p}{2\lambda} 
        + \frac{1}{2} \left(\eta_\infty 
          + \frac{\theta_0}{\lambda} 
          - \frac{\eta_0 t}{\lambda^2}\right). 
    \label{eq:PIII-mu-pq} 
\eeq
This relation, again, can be used to define 
a canonical transformation: 

\begin{theorem}
(\ref{eq:PIII-lam-q}) and (\ref{eq:PIII-mu-pq}) 
define a canonical transformation that connects 
$\PIII$ and the $\PIII$-version of Manin's 
Hamiltonian system.  The canonical coordinates 
and the Hamiltonians of the two systems obey 
the equation 
\beq
    \mu d\lambda - Hdt 
    = \frac{1}{2} \left(pdq - \calH\frac{dt}{t}\right) 
    + \mbox{\rm exact form}. 
\eeq
\end{theorem}

\proof
Since $d\lambda$ and $dq$ are connected by 
the relation 
\beq
    d\lambda = \lambda dq, 
    \nonumber 
\eeq
$\mu d\lambda$ can be written 
\beq
    \mu d\lambda 
    &=& \frac{1}{2}pdq 
      + \frac{1}{2}\left(\eta_\infty 
          + \frac{\theta_0}{\lambda} 
          - \frac{\eta_0 t}{\lambda^2}\right)d\lambda 
    \nonumber \\
    &=& \frac{1}{2}pdq 
      - \frac{\eta_0}{2\lambda}dt 
      + \frac{1}{2}d\left(\eta_\infty\lambda 
          + \theta_0\log\lambda 
          + \frac{\eta_0 t}{\lambda}\right), 
    \nonumber 
\eeq
so that 
\beq
    \mu d\lambda - Hdt 
    = \frac{1}{2}\left(pdq 
        - \tilde{\calH}\frac{dt}{t}\right) 
    + \mbox{\rm exact form}, 
\eeq
where 
\beq
    \tilde{\calH} = 2Ht + \frac{\eta_0 t}{\lambda}. 
\eeq
We can convert the transformed Hamiltonian 
$\tilde{\calH}$ to a normal form as 
\beq
    \tilde{\calH} 
    &=& 2\lambda^2\left[
          \mu - \frac{1}{2}\left( 
            \eta_\infty + \frac{\eta_0}{\lambda} 
            - \frac{\eta_0 t}{\lambda^2}\right)
        \right]^2 
    \nonumber \\
    && + 2\lambda^2\left[ 
          - \frac{1}{2}\left( 
             \eta_\infty + \frac{\eta_0}{\lambda} 
             - \frac{\eta_0 t}{\lambda^2}\right)^2 
          + \frac{\eta_\infty(\theta_0 + \theta_\infty)}
            {2\lambda} \right] 
       + \frac{\eta_0 t}{\lambda} 
    \nonumber \\
    &=& \frac{p^2}{2} + \tilde{V}(q), 
\eeq
where 
\beq
    \tilde{V}(q) 
    &=& - \frac{\lambda^2}{2}\left( 
          \eta_\infty + \frac{\theta_0}{\lambda} 
          - \frac{\eta_0 t}{\lambda^2} \right)^2 
       + \eta_\infty(\theta_0 + \theta_\infty)\lambda 
       + \frac{\eta_0 t}{\lambda} 
    \nonumber \\ 
    &=&   \eta_\infty\theta_\infty e^q 
        + \eta_0(\theta_0 + 1) t e^{-q} 
        - \frac{\eta_\infty^2}{2} e^{2q} 
        - \frac{\eta_0^2 t^2}{2} e^{-2q} 
    \nonumber \\
    && + \mbox{\rm function of $t$ only}. 
\eeq
Thus, apart from the last irrelevant term, 
$\tilde{V}(q)$ coincides with the potential $V(q)$ 
in the statement of the theorem. 
\qed

\subsection{Status of $\PII$ and $\PI$} 

Let us turn to $\PII$ and $\PI$.  The Hamiltonian 
of $\PI$ is already of the normal form 
$\calH = \frac{p^2}{2} + V(q)$ with  $\lambda = q$, 
$\mu = p$ and $H = \calH$.  Although this is not 
the case for $\PII$, one can directly find a canonical 
transformation that converts the Hamiltonian $H$ to 
a normal form: 

\begin{theorem}
A $\PII$-version of Manin's Hamiltonian system 
is defined by the Hamiltonian 
\beq
    \calH= \frac{p^2}{2} 
      - \frac{1}{2}\left(q^2 + \frac{t}{2}\right)^2 
      - \alpha q. 
\eeq
This system is connected with $\PII$ by the 
canonical transformation 
\beq
    \lambda = q, \quad 
    \mu = p + \lambda^2 + \frac{t}{2}. 
\eeq
The canonical coordinates and the Hamiltonians 
of the two systems obey the equation 
\beq
    \mu d\lambda - Hdt 
    = pdq - \calH dt 
    + \mbox{\rm exact form}. 
\eeq
\end{theorem}

\proof 
The foregoing relation between $(\lambda,\mu)$ 
and $(q,p)$ implies that 
\beq
    \mu d\lambda 
    = pdq + \left(\lambda^2 + \frac{t}{2}\right)d\lambda 
    = pdq - \frac{\lambda}{2}dt 
      + d\left(\frac{\lambda^3}{3} 
         + \frac{t\lambda}{2}\right), 
    \nonumber 
\eeq
so that 
\beq
    \mu d\lambda - Hdt 
    = pdq - \tilde{\calH}dt + \mbox{\rm exact form}, 
\eeq
where 
\beq
   \tilde{\calH} 
   &=& H + \frac{\lambda}{2} 
   \nonumber \\
   &=& \frac{1}{2}\left[\mu - \left(\lambda^2 + 
                  \frac{t}{2}\right)\right]^2 
     - \frac{1}{2}\left(\lambda^2 + \frac{t}{2}\right)^2 
     - \left(\alpha + \frac{1}{2}\right) \lambda 
     + \frac{\lambda}{2} 
   \nonumber \\
   &=& \frac{p^2}{2} 
     - \frac{1}{2}\left(q^2 + \frac{t}{2}\right)^2 
     - \alpha q. 
\eeq
This is nothing but the Hamiltonian in the statement 
of the theorem. 
\qed

\section{Multi-component Painlev\'e equations}

\subsection{Inozemtsev Hamiltonians of higher rank}

The rank $\ell$ version of Inozemtsev's Hamiltonians 
have $\ell$ coordinates $q_1,\ldots,q_\ell$ and 
canonical conjugate momenta $p_1,\ldots,p_\ell$.  
The Hamiltonians of the elliptic, hyperbolic 
and rational models take the following form 
\cite{bib:In-Me,bib:Inozemtsev,bib:Le-Wo}: 
\begin{itemize}
\item Elliptic model:
\beq
    \calH 
    = \sum_{j=1}^\ell \left(\frac{p_j^2}{2} 
        + \sum_{n=0}^3 g_n^2 \wp(q_j + \omega_n)
        \right)
    + g_4^2 \sum_{j\not= k} \left(
        \wp(q_j - q_k) + \wp(q_j + q_k)\right). 
    \nonumber 
\eeq
\item Hyperbolic model: 
\beq
    \calH 
    &=& \sum_{j=1}^\ell \left(\frac{p_j^2}{2} 
        + \frac{g_0^2}{\sinh^2(q_j/2)} 
        + \frac{g_1^2}{\cosh^2(q_j/2)} 
        + g_2^2 \cosh(q_j) 
        + g_3^2 \cosh(2q_j) \right) 
    \nonumber \\
    && + g_4^2 \sum_{j\not= k} \left(
         \frac{1}{\sinh^2((q_j - q_k)/2)} 
         + \frac{1}{\sinh^2((q_j + q_k)/2)}\right). 
    \nonumber 
\eeq
\item Rational model: 
\beq
    \calH 
    = \sum_{j=1}^\ell \left(\frac{p_j^2}{2} 
        + g_0^2 q_j^6 + g_1^2 q_j^4 
        + g_2^2 q_j^2 + g_3^2 q_j^{-2} \right) 
      + g_4^2 \sum_{j\not= k} \left( 
         \frac{1}{(q_j - q_k)^2} 
         + \frac{1}{(q_j + q_k)^2} \right). 
    \nonumber 
\eeq
\end{itemize}
Here $g_0,g_1,g_2,g_3$ and $g_4$ are coupling 
constants.  The Painlev\'e-Calogero correspondence 
for $\PIII$, $\PII$ and $\PI$ suggests the existence 
of further degeneration of these models.  

The goal of this section is to extend the 
the Painlev\'e-Calogero correspondence to 
these higher rank models.  Since a complete 
exposition will become inevitably lengthy, 
we shall illustrate the elliptic and 
hyperbolic models in detail, leaving 
the other cases rather sketchy.  
The strategy is as follows:  The point of 
departure is the Hamiltonian of Inozemtsev's 
rank $\ell$ elliptic model.  This gives rise 
to a rank $\ell$ version of Manin's equation.  
Starting with this non-autonomous Hamiltonian 
system, we seek for an analogue of the 
degeneration process for the Painlev\'e 
equations. We can thus obtain six types of 
non-autonomous Hamiltonian systems.  At each 
stage of the degeneration process, we confirm 
that the non-autonomous Hamiltonian system 
on the Calogero side can be mapped, by a canonical 
transformation, to a multi-component analogue of 
the Painlev\'e equation of the corresponding type.

\subsection{Elliptic model and multi-component $\PVI$} 

We now consider the non-autonomous Hamiltonian system 
\beq
    2\pi i\frac{dq_j}{d\tau} = \frac{\rd\calH}{\rd p_j}, 
    \quad 
    2\pi i\frac{dp_j}{d\tau} = - \frac{\rd\calH}{\rd q_j} 
    \label{eq:multi-PVI-Manin}
\eeq
defined by the Hamiltonian of Inozemtsev's elliptic 
model.  This is a rank $\ell$ version of Manin's 
equation. This non-autonomous system is known to 
describe a family of isomonodromic deformations 
on the torus \cite{bib:Ta-ecm}.  

An honest generalization of the canonical 
transformation for the case of $\ell = 1$ leads 
to a multi-component version of $\PVI$ as follows: 

\begin{theorem}
The time-dependent canonical transformation defined by 
\beq
    \lambda_j &=& \frac{\wp(q_j) - e_1}{e_2 - e_1}, 
    \nonumber \\
    \mu_j &=& \frac{e_2 - e_1}{\wp'(q)} p_j 
    + \frac{2\pi i(e_2 - e_1)^2}{\wp'(q_j)^2} f_\tau(q_j) 
    \nonumber \\
    && + \frac{e_2 - e_1}{2} \left( 
         \frac{\kappa_0}{\wp(q_j) - e_1} 
         + \frac{\kappa_1}{\wp(q_j) - e_2} 
         + \frac{\theta - 1}{\wp(q_j) - e_3} \right), 
\eeq
and 
\beq 
    t &=& \frac{e_3 - e_1}{e_2 - e_1}. 
\eeq
maps (\ref{eq:multi-PVI-Manin}) to the Hamiltonian system 
\beq
    \frac{d\lambda_j}{dt} = \frac{\rd H}{\rd\mu_j}, 
    \quad 
    \frac{d\mu_j}{dt} = - \frac{\rd H}{\rd\lambda_j} 
\eeq
with the Hamiltonian 
\beq
    H &=& \sum_{j=1}^\ell 
      \frac{\lambda_j(\lambda_j-1)(\lambda_j-t)}{t(t-1)} 
      \left[ 
        \mu_j^2 
        - \left(\frac{\kappa_0}{\lambda_j} 
          + \frac{\kappa_1}{\lambda_j-1}
          + \frac{\theta-1}{\lambda_j-t}
          \right) \mu_j 
        + \frac{\kappa}{\lambda_j(\lambda_j-1)} 
      \right] 
    \nonumber \\
    && + \frac{g_4^2}{2t(t-1)} \sum_{j\not= k}
         \left[
           \frac{\lambda_j(\lambda_j-1)(\lambda_j-t) 
             + \lambda_k(\lambda_k-1)(\lambda_k-t)}
             {8(\lambda_j - \lambda_k)^2} 
           - 2(\lambda_j + \lambda_k) 
         \right]. 
\eeq
\end{theorem}

\proof
The method of proof for the case of $\ell = 1$ 
can be applied to the present case as well, 
yielding the equality 
\beq
    \sum_{j=1}^\ell p_j dq_j 
    - \calH \frac{d\tau}{2\pi i} 
    = \sum_{j=1}^\ell \mu_j d\lambda_j 
    - \tilde{H} dt 
    + \mbox{\rm exact form}, 
\eeq
where 
\beq
    \tilde{H} 
    &=& \sum_{j=1}^\ell 
        \frac{\lambda_j(\lambda_j-1)(\lambda_j-t)}
             {t(t-1)} 
        \left[ 
            \mu_j^2 
            - \left(\frac{\kappa_0}{\lambda_j} 
              + \frac{\kappa_1}{\lambda_j-1} 
              + \frac{\theta-1}{\lambda_j-t} 
              \right) \mu_j 
            + \frac{\kappa}{\lambda_j(\lambda_j-1)} 
        \right] 
    \nonumber \\
    && + \frac{g_4^2}{2t(t-1)(e_2-e_1)} 
         \sum_{j\not= k} \Bigl( \wp(q_j-q_k) 
           + \wp(q_j+q_k) \Bigr). 
\eeq
What remains is to express the ``two-body potential'' 
part  in terms of $\lambda_j$.  To this end, 
let us recall the addition formula 
\beq
    \wp(u - v) + \wp(u + v) 
    = - 2\wp(u) - 2\wp(v) 
    + \frac{\wp'(u)^2 + \wp'(v)^2}{2(\wp(u) - \wp(v))^2} 
\eeq
of the $\wp$-function.  Applying it to the case 
where $(u,v) = (\lambda_j,\lambda_k)$, and 
substituting 
\beq
    \wp(q_j) &=& e_1 + (e_2 - e_1)\lambda_j, 
    \nonumber \\
    \wp(q_k) &=& e_1 + (e_2 - e_1)\lambda_k, 
    \nonumber \\
    \wp'(q_j)^2 &=& \frac{(e_2-e_1)^3}{4}
                \lambda_j(\lambda_j-1)(\lambda_j-t), 
    \nonumber \\
    \wp'(q_k)^2 &=& \frac{(e_2-e_1)^3}{4} 
                \lambda_k(\lambda_k-1)(\lambda_k-t), 
    \nonumber 
\eeq
we can rewrite the two-body potential terms as 
\beq
    \wp(q_j - q_k) + \wp(q_j + q_k) 
    &=& - 2 (e_1 + (e_2-e_1)\lambda_j) 
        - 2 (e_1 + (e_2-e_1)\lambda_k) 
    \nonumber \\
    && + \frac{(e_2-e_1)^3}{8} \cdot 
         \frac{\lambda_j(\lambda_j-1)(\lambda_j-t) 
             + \lambda_k(\lambda_k-1)(\lambda_k-t)} 
         {(e_1 + (e_2-e_1)\lambda_j 
                - e_1 - (e_2-e_1)\lambda_k)^2}
    \nonumber \\
    &=& - 4e_1 - 2(e_2-e_1)(\lambda_j+\lambda_k) 
    \nonumber \\
    && + \frac{e_2-e_1}{8} \cdot 
         \frac{\lambda_j(\lambda_j-1)(\lambda_j-t) 
             + \lambda_k(\lambda_k-1)(\lambda_k-t)} 
              {(\lambda_j-\lambda_k)^2} . 
\eeq 
The first term $-4e_1$ is non-dynamical, 
thereby negligible (i.e., can be absorbed by 
the ``exact form'' part).  Removing these terms 
from $\tilde{H}$, we obtain the Hamiltonian $H$. 
\qed

\subsection{Degeneration of elliptic model to 
hyperbolic model}

The degeneration of the elliptic model is 
achieved by letting $\Im\tau \to +\infty$.  
Like the degeneration process from $\PVI$ 
to $\PV$, this is a kind of scaling limit, 
namely, the coupling constants $g_n$ and 
the elliptic modulus $\tau$ have to be suitably 
rescaled.  To this end, we have to understand 
the asymptotic behavior of the constants 
$e_1,e_2,e_3$ and the $\wp$-function in the 
limit as $\Im\tau \to +\infty$.  All necessary 
data are collected in Appendix B. For instance, 
the asymptotic expression of $e_1,e_2$ and $e_3$ 
imply that 

\beq
    t = 1 + \frac{e_3 - e_2}{e_2 - e_1} 
      = 1 + 16\pi^2 e^{\pi i\tau} + O(e^{2\pi i\tau}). 
\eeq
This is indeed consistent with the scaling 
rule $t = 1 + \epsilon \tilde{t}$ in the 
degeneration process of $\PVI$ to $\PV$.  

Having these data, we now rescale the coupling 
constants and the elliptic modulus as 
\beq
    g_0^2 = \tilde{g}_0^2, \quad 
    g_1^2 = \tilde{g}_1^2, \quad 
    g_2^2 = \frac{\tilde{g}_2^2}{\epsilon} 
          + \frac{\tilde{g}_3^2}{\epsilon^2}, \quad 
    g_3^3 = \frac{\tilde{g}_3^2}{\epsilon^2}, \quad 
    g_4^2 = \tilde{g}_4^2 
\eeq
and 
\beq
    16e^{\pi i\tau} = \epsilon \tilde{t}, 
\eeq
and consider the limit as $\epsilon \to 0$ 
while leaving $\tilde{g}_n$ and $\tilde{t}$ 
finite.  Note that letting $\epsilon \to 0$ 
amounts to letting $\Im\tau \to +\infty$.  

The asymptotic expression of $\wp(u)$ and 
$\wp(u + \omega_n)$ in Appendix B show that 
the potential $V(q)$ of the elliptic model 
behaves as 
\beq
    V(q) &=& 
    \sum_{j=1}^\ell \left( 
        \frac{\tilde{g}_0^2 \pi^2}{\sin^2(\pi q_j)} 
      + \frac{\tilde{g}_1^2 \pi^2}{\cos^2(\pi q_j)} 
      + \frac{\tilde{g}_2^2 \pi^2 \tilde{t}}{2} 
        \cos(2\pi q_j) 
      - \frac{\tilde{g}_3^2 \pi^2 \tilde{t}^2}{8} 
        \cos(4\pi q_j) 
    \right) 
    \nonumber \\
    && + \tilde{g}_4^2 \sum_{j\not= k} \left( 
            \frac{1}{\sin^2(\pi(q_j - q_k))} 
          + \frac{1}{\sin^2(\pi(q_j + q_k))}
    \right) 
    \nonumber \\
    && + \mbox{\rm function of $\epsilon$ and 
         $\tilde{t}$ only} 
       + O(\epsilon). 
    \nonumber 
\eeq
Thus, removing negligible terms, we obtain the 
following Hamiltonian in the limit: 
\beq
    \tilde{\calH} &=& 
    \sum_{j=1}^\ell \left( 
        \frac{p_j^2}{2} 
      + \frac{\tilde{g}_0^2 \pi^2}{\sin^2(\pi q_j)} 
      + \frac{\tilde{g}_1^2 \pi^2}{\cos^2(\pi q_j)} 
      + \frac{\tilde{g}_2^2 \pi^2 \tilde{t}}{2} 
        \cos(2\pi q_j) 
      - \frac{\tilde{g}_3^2 \pi^2 \tilde{t}^2}{8} 
        \cos(4\pi q_j) 
    \right) 
    \nonumber \\
    && 
    + \tilde{g}_4^2 \sum_{j\not= k} \left( 
            \frac{1}{\sin^2(\pi(q_j - q_k))} 
          + \frac{1}{\sin^2(\pi(q_j + q_k))}
      \right). 
\eeq

The asymptotic expression of $t$ determines the 
equation of motion in the limit. In fact, since 
\beq
    \frac{d\tau}{dt} 
    = \frac{\pi}{t(t - 1)(e_2 - e_1)} 
    = \frac{\pi i}{(1 + \epsilon\tilde{t})
        (-\epsilon\tilde{t}) (-\pi^2 + O(\epsilon))} 
    \nonumber 
\eeq
and 
\beq
    2\pi i\frac{d}{d\tau} 
    = 2\pi i \frac{dt}{d\tau} \frac{d\tilde{t}}{dt} 
      \frac{d}{dt} 
    = \Bigl(2\pi^2\tilde{t} + O(\epsilon^2)\Bigr) 
      \frac{d}{dt}, 
    \nonumber 
\eeq
we find that the equations of motion take the 
following form: 
\beq
    2\pi^2\tilde{t}\frac{dq_j}{d\tilde{t}} 
    = \frac{\rd\tilde{\calH}}{\rd p_j},  \quad 
    2\pi^2\tilde{t}\frac{dp_j}{d\tilde{t}} 
    = - \frac{\rd\tilde{\calH}}{\rd q_j}.  
\eeq

The final step is to rescale the variables 
and the Hamiltonian as  
\beq
   q_j \to \frac{q_j}{2\pi i}, \quad 
   p_j \to \pi iq_j, \quad 
   \tilde{\calH} \to - \pi^2\tilde{\calH}, 
\eeq
and to rename $\tilde{t}$ and $\tilde{\calH}$ 
to $t$ and $\calH$.  Let us also define the 
new constants 
\beq
    \alpha = - \frac{\tilde{g}_0^2}{2}, \quad 
    \beta = \frac{\tilde{g}_1^2}{2}, \quad 
    \gamma = - \frac{\tilde{g}_2^2}{2}, \quad 
    \delta = \frac{\tilde{g}_3^2}{2}, 
\eeq
which are to be identified with the four 
parameters of $\PV$. The outcome is the 
non-autonomous Hamiltonian system 
\beq
    t \frac{dq_j}{dt} = \frac{\rd\calH}{\rd p_j}, \quad 
    t \frac{dp_j}{dt} = - \frac{\rd\calH}{\rd q_j}
    \label{eq:multi-PV-Manin}
\eeq
with the Hamiltonian 
\beq
    \calH &=& 
    \sum_{j=1}^\ell \left( 
        \frac{p_j^2}{2} 
      - \frac{\alpha}{\sinh^2(q_j/2)} 
      - \frac{\beta}{\cosh^2(q_j/2)} 
      + \frac{\gamma t}{2} \cosh(q_j) 
      + \frac{\delta t^2}{8} \cosh(2q_j) 
    \right) 
    \nonumber \\
    && 
    + g_4^2 \sum_{j\not= k} \left( 
          \frac{1}{\sinh^2((q_j - q_k)/2)} 
        + \frac{1}{\sinh^2((q_j + q_k)/2)} 
      \right). 
\eeq
This gives a rank $\ell$ version of the 
non-autonomous Hamiltonian system on the 
Calogero side of $\PV$.  Note that the 
Hamiltoian is essentially the same as 
the Hamiltonian of Inozemtsev's hyperbolic 
model, except that the effective coupling 
constants are now time-dependent.  

\remark
The foregoing prescription of scaling limit 
of the coupling constants and the elliptic 
modulus is reminiscent of ``renormalization'' 
in quantum field theories.  In this analogy, 
one can interpret the equations of motion 
of the Hamiltonian system as ``renormalization 
group equations'', in which $\tilde{t}$ plays 
the role of a ``mass scale'' parameter.

\subsection{Canonical transformation to 
multi-component $\PV$} 

Again, an honest generalization of the 
canonical transformation for the case 
of $\ell = 1$ leads to a multi-component 
version of $\PV$: 

\begin{theorem}
The time-dependent canonical transformation 
defined by 
\beq
    \sqrt{\lambda_j} &=& - \coth(q_j/2), 
    \nonumber \\
    \mu_j 
    &=& \frac{p_j}{2\sqrt{\lambda_j}(\lambda_j - 1)} 
      + \frac{1}{2}\left(\frac{\kappa_0}{\lambda_j} 
        + \frac{\theta_1}{\lambda_j - 1} 
        - \frac{\eta_1 t}{(\lambda_j -1)^2} \right) 
\eeq
maps (\ref{eq:multi-PV-Manin}) to the Hamiltonian 
system 
\beq 
    \frac{d\lambda_j}{dt} = \frac{\rd H}{\rd \mu_j}, \quad 
    \frac{d\mu_j}{dt} = - \frac{\rd H}{\rd \lambda_j} 
\eeq
with the Hamiltonian 
\beq
    H &=& 
      \sum_{j=1}^\ell 
      \frac{\lambda_j(\lambda_j - 1)^2}{t} 
      \left[ 
        \mu_j^2 
        - \left(\frac{\kappa_0}{\lambda_j} 
            + \frac{\theta_1}{\lambda_j - 1} 
            - \frac{\eta_1 t}{(\lambda_j - 1)^2} 
          \right) \mu_j 
        + \frac{\kappa}{\lambda_j(\lambda_j - 1)} 
      \right]   
    \nonumber \\
    && + \frac{g_4^2}{2t} 
         \sum_{j\not= k} 
         \frac{2(\lambda_j-1)(\lambda_k-1)(\lambda_j+\lambda_k)} 
              {(\lambda_j-\lambda_k)^2}. 
\eeq
\end{theorem}

\proof 
The method of proof for the case of $\ell = 1$ 
can be used as it is.  The outcome is the equality 
\beq
    \sum_{j=1}^\ell p_j dq_j - \calH \frac{dt}{t} 
    = 2\left(\sum_{j=1}^\ell \mu_j d\lambda_j 
        - H dt \right) 
    + \mbox{\rm exact form}, 
\eeq
where 
\beq
    H &=& 
      \sum_{j=1}^\ell 
      \frac{\lambda_j(\lambda_j - 1)^2}{t} 
      \left[ 
        \mu_j^2 
        - \left(\frac{\kappa_0}{\lambda_j} 
            + \frac{\theta_1}{\lambda_j - 1} 
            - \frac{\eta_1 t}{(\lambda_j - 1)^2} 
          \right) \mu_j 
        + \frac{\kappa}{\lambda_j(\lambda_j - 1)} 
      \right]   
    \nonumber \\
    && + \frac{g_4^2}{2t} \sum_{j\not= k} 
         \left( \frac{1}{\sinh^2((q_j-q_k)/2)} 
         + \frac{1}{\sinh^2((q_j+q_k)/2)} \right). 
\eeq
The two-body potential part can be rewritten 
by use of the identity 
\beq
    \frac{1}{\sinh^2(u-v)} + \frac{1}{\sinh^2(u+v)} 
    = 4 \frac{\cosh(2u)\cosh(2v) - 1} 
             {(\cosh(2u) - \cosh(2v))^2}. 
\eeq
Substituting $u = q_j/2$, $v = q_k/2$, and also using 
the equality $\cosh(q_j) = (\lambda_j + 1)/(\lambda_j - 1)$, 
we find that 
\beq
      \frac{1}{\sinh^2((q_j-q_k)/2)} 
      + \frac{1}{\sinh^2((q_j+q_k)/2)} 
    = \frac{2(\lambda_j-1)(\lambda_k-1)(\lambda_j+\lambda_k)} 
           {(\lambda_j - \lambda_k)^2}, 
\eeq
which gives the two-body potential term in $H$. 
\qed

\subsection{Other models}

The degeneration process can be further continued, 
and leads to four more models that correspond to 
a multi-component version of $\PIV$, $\PIII$, 
$\PII$ and $\PI$.  Since  the details of derivation 
are more or less parallel, we show the final results 
only.  The Hamiltonian of each model, like those in 
the foregoing cases, becomes a sum of $\ell$ copies 
of the one-component Hamiltonian and Calogero-like 
two-body potential terms.

\subsubsection{Rational model and multi-component $\PIV$} 

This model can be derived from the hyperbolic model 
by degeneration.  The degeneration process consists 
of putting the variables and the parameters as 
\beq 
    t = 1 + 2\epsilon\tilde{t}, 
    \quad 
    q_j = \pi i + \epsilon^{1/2}\tilde{q_j}, \quad 
    p_j = \frac{\tilde{p_j}}{2\epsilon^{1/2}}, 
\eeq
and 
\beq
    \alpha = \frac{1}{8\epsilon^4}, \quad 
    \beta = \frac{\tilde{\beta}}{4}, \quad 
    \gamma = \frac{1}{4\epsilon^4}, \quad 
    \delta = - \frac{1}{8\epsilon^4} 
           + \frac{\tilde{\alpha}}{2\epsilon^2}, 
\eeq
and letting $\epsilon \to 0$ while leaving the 
``renormalized'' quantities $\tilde{t}$, etc. 
finite. 

The equations of motion of this model takes 
the canonical form 
\beq
    \frac{dq_j}{dt} = \frac{\rd\calH}{\rd p_j}, \quad 
    \frac{dp_j}{dt} = - \frac{\rd\calH}{\rd q_j} 
\eeq
with the Hamiltonian 
\beq
    \calH 
    &=& \sum_{j=1}^\ell 
        \left[ 
        \frac{p_j^2}{2} 
        - \frac{1}{2} \left(\frac{q_j}{2}\right)^6 
        - 2t \left(\frac{q_j}{2}\right)^4 
        - 2(t^2 - \alpha)\left(\frac{q_j}{2}\right)^2 
        + \beta \left(\frac{q_j}{2}\right)^{-2} 
        \right]
    \nonumber \\
    && + g_4^2 \sum_{j\not= k} \left( 
           \frac{1}{(q_j - q_k)^2}  
           + \frac{1}{(q_j + q_k)^2} \right). 
\eeq
The canonical transformation defined by 
\beq
    \lambda_j = \left(\frac{q_j}{2}\right)^2, \quad 
    \mu_j = \frac{p_j}{4\sqrt{\lambda_j}} 
          + \frac{1}{4} \left(\lambda_j + 2t 
            + \frac{2\kappa_0}{\lambda_j}\right) 
\eeq
maps the foregoing non-autonomous system to 
the Hamiltonian system 
\beq
    \frac{d\lambda_j}{dt} = \frac{\rd H}{\rd\mu_j}, 
    \quad 
    \frac{d\mu_j}{dt} = - \frac{\rd H}{\rd\lambda_j} 
\eeq
with the Hamiltonian 
\beq
    H = \sum_{j=1}^\ell 
          2\lambda_j^2 \left[ 
            \mu_j^2 
            - \left(\frac{\lambda_j}{2} + t 
              + \frac{\kappa_0}{\lambda}\right) \mu_j 
            + \frac{\theta_0}{2} 
          \right] 
       + \frac{g_4^2}{4}
         \sum_{j\not= k} 
           \frac{2(\lambda_j + \lambda_k)}
                {(\lambda_j - \lambda_k)^2}.  
\eeq

\subsubsection{Exponential-hyperbolic model and 
multi-component $\PIII$} 

This model, too, can be derived from the hyperbolic 
model by degeneration.   This degeneration is 
achieved by the putting the variables and the 
parameters as 
\beq
    q_j = - \tilde{q_j} - \log \frac{\epsilon}{4}, \quad 
    p_j = - \tilde{p_j}, 
\eeq
and 
\beq 
    \alpha = \frac{\tilde{\alpha}}{4\epsilon} 
           + \frac{\tilde{\gamma}}{8\epsilon^2}, \quad 
    \beta = - \frac{\tilde{\gamma}}{8\epsilon^2}, \quad 
    \gamma = \frac{\tilde{\beta}\epsilon}{4}, \quad 
    \delta = \frac{\tilde{\delta}\epsilon^2}{8}, 
\eeq
and letting $\epsilon \to 0$.

The equations of motion of this model takes the 
canonical form 
\beq
    t \frac{dq_j}{dt} = \frac{\rd\calH}{\rd p_j}, 
    \quad 
    t \frac{dp_j}{dt} = - \frac{\rd\calH}{\rd q_j} 
\eeq
with the Hamiltonian 
\beq
    \calH 
    &=& \sum_{j=1}^\ell 
        \left( \frac{p_j^2}{2} 
          - \frac{\alpha}{4} e^{q_j} 
          + \frac{\beta t}{4} e^{-q_j} 
          - \frac{\gamma}{8} e^{2q_j} 
          + \frac{\delta t^2}{8} e^{-2q_j} \right) 
    \nonumber \\
    && + g_4^2 \sum_{j\not= k} 
         \frac{1}{\sinh^2((q_j - q_k)/2)}. 
\eeq
The canonical transformation defined by 
\beq
   \lambda_j = e^{q_j}, \quad 
   \mu_j = \frac{p_j}{2\lambda_j} 
         + \frac{1}{2}\left( \eta_\infty 
           + \frac{\theta_0}{\lambda_j} 
           - \frac{\eta_0 t}{\lambda_j^2} \right) 
\eeq
maps the foregoing onn-autonomous system to 
the Hamiltonian system 
\beq
    \frac{d\lambda_j}{dt} = \frac{\rd H}{\rd \mu_j}, 
    \quad 
    \frac{d\mu_j}{dt} = - \frac{\rd H}{\rd \lambda_j} 
\eeq
with the Hamiltonian 
\beq
    H = \sum_{j=1}^\ell 
          \frac{\lambda_j^2}{t} \left[ 
            \mu_j^2 
            - \left(\eta_\infty + \frac{\theta_0}{\lambda_j} 
              - \frac{\eta_0 t}{\lambda_j^2} \right) \mu_j 
            + \frac{\eta_\infty(\theta_0 + \theta_\infty)} 
                   {2\lambda_j} 
          \right] 
     + \frac{g_4^2}{2t} \sum_{j\not= k} 
         \frac{4\lambda_j\lambda_k}
              {(\lambda_j - \lambda_k)^2}.  
\eeq

\subsubsection{Second rational model and 
multi-component $\PII$} 

This model can be derived from {\it both} the rational 
model and the exponential-hyperbolic model by degeneration. 
For the degeneration from the rational model, we 
write the variables and the parameters as 
\beq
    t = \frac{-1 + 4^{-1/3}\epsilon^4\tilde{t}}{\epsilon}, 
    \quad 
    \frac{q_j}{2} 
    = \frac{1 + 2^{-1/3}\epsilon^2\tilde{q_j}}{\epsilon^{3/2}}, 
    \quad 
    p_j = \frac{4^{2/3}\tilde{p}_j}{\epsilon^{1/2}}  
\eeq
and 
\beq
    \alpha = - 2\tilde{\alpha} - \frac{1}{2\epsilon^6}, 
    \quad 
    \beta = - \frac{1}{2\epsilon^{12}}, 
\eeq
and let $\epsilon \to 0$.  The degeneration from 
the exponential-hyperbolic model is similarly achieved 
by putting 
\beq
    t = 1 + 2 \epsilon^2 \tilde{t}, \quad 
    q_j = 2 \epsilon \tilde{q_j}, \quad 
    p_j = \frac{\tilde{p_j}}{\epsilon}, 
\eeq
and 
\beq
    \alpha = - \frac{1}{2\epsilon^6}, \quad 
    \beta = \frac{1 + 4\epsilon^3\tilde{\alpha}}{2\epsilon^6}, 
    \quad 
    \gamma = \frac{1}{4\epsilon^6}, \quad 
    \delta = - \frac{1}{4\epsilon^6}, 
\eeq
and again letting $\epsilon \to 0$.  

The equations of motion of this model takes the 
canonical form 
\beq
    \frac{dq_j}{dt} = \frac{\rd\calH}{\rd p_j}, \quad 
    \frac{dp_j}{dt} = - \frac{\rd\calH}{\rd q_j} 
\eeq
with the Hamiltonian 
\beq
    \calH 
    = \sum_{j=1}^\ell \left[ 
          \frac{p_j^2}{2} 
          - \frac{1}{2}\left(q_j^2 + \frac{t}{2}\right)^2 
          - \alpha q_j \right] 
    + g_4^2 \sum_{j\not= k} 
         \frac{1}{(q_j - q_k)^2}. 
\eeq
The canonical transformation defined by 
\beq
    \lambda_j = q_j, \quad 
    \mu_j = p_j + \lambda_j^2 + \frac{t}{2} 
\eeq
maps the foregoing non-autonomous system to the 
Hamiltonian system 
\beq
    \frac{d\lambda_j}{dt} = \frac{\rd H}{\rd\mu_j}, 
    \quad 
    \frac{d\mu_j}{dt} = - \frac{\rd H}{\rd\lambda_j} 
\eeq
with the Hamiltonian 
\beq
    H =  \sum_{j=1}^\ell \left[ 
          \frac{\mu_j^2}{2} 
          - \left(\lambda_j^2 + \frac{t}{2}\right)\mu_j 
          - \left(\alpha + \frac{1}{2}\right) \lambda_j 
         \right] 
       + g_4^2 \sum_{j\not= k} 
         \frac{1}{(\lambda_j - \lambda_k)^2}. 
\eeq

\subsubsection{Multi-component $\PI$} 

This model can be derived from the second rational 
model, and takes the {\it same} form on both the 
Painlev\'e and Calogero sides.  The degeneration process 
is achieved by putting 
\beq
    t = \frac{-6 + \epsilon^{12}\tilde{t}}{\epsilon^{10}}, 
    \quad 
    q_j = \frac{1 + \epsilon^6\tilde{q_j}}{\epsilon^5}, 
    \quad 
    p_j = \frac{\tilde{p_j}}{\epsilon}, 
    \quad 
    \alpha = {4}{\epsilon^{15}} 
\eeq
and letting $\epsilon \to 0$.  The equations of motion 
takes the canonical form 
\beq
    \frac{dq_j}{dt} = \frac{\rd H}{\rd p_j}, \quad 
    \frac{dp_j}{dt} = - \frac{\rd H}{\rd q_j} 
\eeq
with the Hamiltonian 
\beq
    H = \sum_{j=1}^\ell 
        \left( \frac{p_j^2}{2} - 2q_j^3 - tq_j \right) 
      + g_4^2 \sum_{j\not= k} \frac{1}{(q_j - q_k)^2}. 
\eeq

\section{Concluding remarks}

We have shown that the Painlev\'e-Calogero 
correspondence persists for all the six 
Painlev\'e equations and their multi-component 
generalizations. The Calogero side of this 
correspondence is a non-autonomous version of 
Inozemtsev's elliptic model and its various 
degenerations. Those for $\PV$ and $\PIV$ are 
a non-autonomous version of Inozemtsev's 
hyperbolic and rational models. The others 
corresponding to $\PIII$, $\PII$ and $\PI$ are 
further degenerations of the hyperbolic and 
rational models.  The pattern of degeneration 
on the Calogero side repeats the degeneration 
diagram 
\beq
  \begin{array}{ccccccc}
    \PVI & \longrightarrow & \PV & \longrightarrow & 
    \PIV & & \\
         &                 & \downarrow & & 
    \downarrow & & \\
         &                 & \PIII & \longrightarrow & 
    \PII & \longrightarrow & \PI 
  \end{array}
  \nonumber 
\eeq
of the Painlev\'e equations.  

This picture applies to the autonomous systems 
as well.  Actually, such degeneration relations 
in the autonomous case have been more or less well 
known to experts of Calogero-Moser systems (see the 
Introduction of van~Diejen's paper \cite{bib:vanDiejen}). 
The autonomous systems are defined by a Hamiltonian 
of the same form with the time-dependent coupling 
constants being replaced by absolute constants (except 
for the elliptic model, in which case an independent 
time variable is introduced).  Those in the position 
of the first row of the degeneration diagram are, 
of course, Inozemtsev's elliptic, hyperbolic and 
rational models (see Section 5).  Those in the 
position of $\PIII$ and $\PII$ are defined by the 
following Hamiltonians: 
\begin{itemize}
\item Exponential-hyperbolic model: 
\beq
    \calH 
    = \sum_{j=1}^\ell \left(
       \frac{p_j^2}{2} 
       + g_0^2 e^{q_j} + g_1^2 e^{2q_j} 
       + g_2^2 e^{-q_j} + g_3^2 e^{-2q_j} \right) 
    + g_4^2 \sum_{j\not= k} 
       \frac{1}{\sinh^2((q_j - q_k)/2)}. 
    \nonumber 
\eeq
\item Second rational model: 
\beq
    \calH 
    = \sum_{j=1}^\ell \left(
       \frac{p_j^2}{2} 
       + g_0^2 q_j^4 + g_1^2 q_j^3 
       + g_2^2 q_j^2 + g_3^2 q_j \right) 
    + g_4^2 \sum_{j\not= k} 
       \frac{1}{(q_j - q_k)^2}. 
    \nonumber 
\eeq
\end{itemize}
The Hamiltonian in the position of $\PI$ is redundant 
in the automonous case, because it is a specialization, 
rather than a degeneration, of the last Hamiltonian.  

Note that the Hamiltonian of the second rational 
model is a {\it quartic} perturbation of the usual 
($A_\ell$ type) rational Calogero Hamiltonian.  
According to recent work of Caseiro, Fran\c{c}oise 
and Sasaki \cite{bib:Ca-Fr-Sa}, such a quartic 
(integrable) perturbation always exists for any 
rational Calogero-Moser system. Inozemtsev's 
rational model, which is a {\it sextic} perturbation 
of the $D_\ell$ type rational Calogero-Moser system, 
might admit a similar interpretation.  

Back to the Painlev\'e equations, the extended 
Painlev\'e-Calogero correspondence raises many 
interesting problems.  A central issue will be 
to find an isomonodromic description of the 
multi-component Painlev\'e equations.  If 
such an isomonodromic description does exist, 
it should be related to a new geometric structure.

\subsection*{Acknwlegements}

I am grateful to Marta Mazzocco, Davide Guzzetti, 
Kazuo Okamoto, Ryu Sasaki, Shun Shimomura and 
Jan Felipe van~Diejen, for useful comments. 
This work was partly supported by the Grant-in-Aid 
for Scientific Researches (No. 10640165) from 
the Ministry of Education, Science and Culture. 

%%%%%%%%%%%%%%%%%%%%
%%%% Appendices %%%%
%%%%%%%%%%%%%%%%%%%%
\appendix 
\renewcommand{\theequation}{\Alph{section}.\arabic{equation}}

\section{Proof of (\ref{eq:term-A})}
\setcounter{equation}{0}

Let us introduce the two auxiliary functions 
\beq
    g(u) = \frac{f_\tau(u)}{f'(u)}, \quad 
    h(u) = \frac{\vartheta'(u + \omega_1)}
                {\vartheta(u + \omega_1)}, 
\eeq
associated with the function 
\beq
    f(u) = \frac{\wp(u) - e_1}{e_2 - e_1} 
\eeq
and the standard elliptic theta function 
\beq
    \vartheta(u) 
    = \sum_{n=-\infty}^\infty 
      \exp(\pi i\tau n^2 + 2\pi inu). 
\eeq

\begin{lemma}
$g(u)$ is a meromorphic function on the 
$u$-plane with additive quasi-periodicity 
\beq
    g(u+1) = g(u), \quad 
    g(u+\tau) = g(u) - 1. 
\eeq
All poles are of the first order and contained 
in the lattice $\omega_3 + \bbZ + \tau\bbZ$.  
Furthermore, $g(u)$ has zeros at $u = 0$ and 
$u = \omega_1$. 
\end{lemma}

\proof
Since $f(u)$ is a doubly periodic function 
with primitive periods 1 and $\tau$, 
$f'(u)$ and $f_\tau(u)$ transform as 
\beq
    f'(u+1) = f'(u), && 
    f'(u+\tau) = f'(u), 
    \nonumber \\
    f_\tau(u+1) = f_\tau(u), && 
    f_\tau(u+\tau) = f_\tau(u) - f'(u) 
    \nonumber 
\eeq
under the shift by $1$ and $\tau$. 
This implies the additive quasi-periodicity 
of $g(u)$. Furthermore, by the construction, 
$g(u)$ is a meromorphic function on the 
$u$-plane, and all possible poles are of 
the first order and located at the points of 
$\omega_k + \bbZ + \tau\bbZ$. Let us examine 
the behavior of $g(u)$ at the representative 
points $u = \omega_0,\omega_1,-\omega_2,\omega_3$: 
\begin{itemize}
\item 
As $u \to \omega_0 = 0$, 
\beq
    f(u) = \frac{1}{(e_2 - e_1)u^2} + O(1), 
    \nonumber 
\eeq
thereby 
\beq
   f'(u) = - \frac{2}{(e_2 - e_1)u^3} + O(1), 
   \quad 
   f_\tau(u) 
   = - \frac{e_{2,\tau} - e_{1,\tau}}
       {(e_2 - e_1)^2 u^2}  + O(1), 
   \nonumber 
\eeq
so that $g(u)$ has rather a zero at $u = 0$: 
\beq
    g(u) = O(u). 
\eeq
\item As $u \to \omega_1 = \frac{1}{2}$, 
\beq
    f(u) 
    &=& \frac{1}{e_2 - e_1}
        \Bigl(\wp(\omega_1) - e_1 
           + \wp'(\omega_1)(u - \omega_1) 
           + O((u - \omega_1)^2) \Bigr) 
    \nonumber \\
    &=& O((u - \omega_1)^2), 
    \nonumber 
\eeq
thereby 
\beq
    f'(u) = O(u -\omega_1), \quad 
    f_\tau(u) = O((u - \omega_1)^2), 
    \nonumber 
\eeq
so that $g(u)$ has another zero at 
$u = \omega_1$: 
\beq
    g(u) = O(u - \omega_1). 
\eeq
\item 
As $u \to - \omega_2 = \frac{1}{2} + \frac{\tau}{2}$, 
\beq
    f(u) 
    &=& \frac{1}{e_2 - e_1}
        \Bigl(\wp(-\omega_2) - e_1 
          + \wp'(-\omega_2)(u + \omega_2) 
          + O((u + \omega_2)^2) \Bigr) 
    \nonumber \\
    &=& O((u + \omega_2)^2), 
    \nonumber 
\eeq
thereby 
\beq
    f'(u) = O(u + \omega_2), \quad 
    f_\tau(u) = O(u + \omega_2), 
    \nonumber 
\eeq
so that $g(u)$ behaves as 
\beq
    g(u) = O(1). 
\eeq
\item 
As $u \to \omega_3 = \frac{\tau}{2}$, 
\beq
    f(u) 
    &=& \frac{1}{e_2 - e_1}
        \Bigl(\wp(\omega_3) - e_1 
          + \wp'(\omega_3)(u - \omega_3) 
          + O((u - \omega_3)^2) \Bigr) 
    \nonumber \\
    &=& t + O((u - \omega_3)^2), 
    \nonumber 
\eeq
thereby 
\beq
    f'(u) = O(u - \omega_3), \quad 
    f_\tau(u) = O(1), 
    \nonumber 
\eeq
so that $g(u)$ turns out to have a pole of 
the first order at $u = \omega_3$: 
\beq
    g(u) = 0((u - \omega_3)^{-1}). 
\eeq
\end{itemize}
The behavior of $g(u)$ at the other 
points of $\omega_n + \bbZ + \tau\bbZ$ 
can be deduced from these results by 
the additive quasi-periodicity of $g(u)$.  
\qed

\begin{lemma}
$h(u)$ is a meromorphic function on the 
$u$-plane with additive quasi-periodicity 
\beq
    h(u + 1) = h(u), \quad 
    h(u + \tau) = h(u) - 2\pi i. 
\eeq
All poles are of the first order and contained 
in the lattice $\omega_3 + \bbZ + \tau\bbZ$.  
Furthermore, $h(u)$ has zeros at $u = 0$ 
and $u = \omega_1$. 
\end{lemma}

\proof 
Let us recall the fundamental properties of  
$\vartheta(u)$: 
\begin{itemize}
\item $\vartheta(u)$ is an entire function on the 
$u$-plane with zeros of the first order at the lattice 
points $\omega_2 + m + n\tau$ ($m,n \in \bbZ$). 
\item $\vartheta(u)$ is quasi-periodic, 
\beq 
    \vartheta(u + 1) = \vartheta(u), \quad 
    \vartheta(u + \tau) 
      = e^{-\pi i \tau - 2\pi i u} \vartheta(u). 
    \nonumber 
\eeq
\item $\theta(u)$ and $\vartheta(u + 1/2)$ are 
even under the reflection $u \to -u$. 
\end{itemize}
All the properties of $h(u)$ in the statement 
of the lemma are an immediate consequence of 
these properties of $\vartheta(u)$. 
\qed

\begin{lemma} 
The function $f(u)$ satisfies the equation 
\beq
      2\pi i\frac{f_\tau(u)}{f'(u)} 
    = \frac{\vartheta'(u + \omega_1)}
      {\vartheta(u + \omega_1)}, 
\eeq
where the prime stands for $\rd/\rd u$. 
\end{lemma}

\proof
The foregoing properties of $g(u)$ and $h(u)$ 
imply the following: 
\begin{itemize}
\item $2\pi i g(u) - h(u)$ is a doubly periodic 
meromorphic function with fundamental period $1$ 
and $\tau$. 
\item All poles of $2\pi ig(u) - h(u)$ are of 
the first order and contained in the lattice 
$\omega_3 + \bbZ + \tau\bbZ$. 
\item $2\pi i g(u) - h(u)$ has zeros at 
$u = 0$ and $u = \omega_1$.  
\end{itemize}
The first two properties imply that 
$2\pi ig(u) - h(u)$ is a constant. 
By the last one, this constant has 
to be zero.  We thus find that 
$2\pi ig(u) - h(u) = 0$. 
\qed

\begin{lemma}
$\vartheta(u)$ satisfies the equation 
\beq
    \Bigl(\log\vartheta(u + \omega_1)\Bigr)'' 
    = - \wp(u + \omega_3) 
    + \mbox{\rm function of $\tau$ only}. 
\eeq
\end{lemma}

\proof
The aforementioned complex analytic properties 
of $\vartheta(u)$ imply the following: 
\begin{itemize}
\item $\Bigl(\log\vartheta(u + \omega_1)\Bigr)''$ 
is a doubly periodic meromorphic function with 
primitive period $1$ and $\tau$. 
\item All poles of this meromorphic function 
are contained in the lattice $\omega_3 + \bbZ 
+ \tau\bbZ$. 
\item As $u \to - \omega_3$, this function 
behaves as 
\beq
    \Bigl(\log\vartheta(u+\omega_1) \Bigr)'' 
    = - \frac{1}{(u + \omega_3)^2} + O(1). 
    \nonumber 
\eeq
\end{itemize}
The function $- \wp(u + \omega_3)$, too, has 
these properties.  Accordingly, their difference 
is a constant function on the $u$-plane, namely, 
a function of $\tau$ only. 
\qed

We now return to the proof of (\ref{eq:term-A}).  
By the third lemma, we have the identity 
\beq
    2\pi i \frac{f_\tau(u)}{f'(u)}du 
    = \frac{\vartheta'(u + \omega_1)}
             {\vartheta(u + \omega_1)} du 
    = \frac{d\vartheta(u + \omega_1)}
             {\vartheta(u + \omega_1)} 
      - \frac{\rd\vartheta(u + \omega_1)/\rd\tau}
             {\vartheta(u + \omega_1)} d\tau 
\eeq
On the other hand, the well known ``heat equation'' 
\beq
    4\pi i\frac{\rd\vartheta(u)}{\rd\tau} 
    = \vartheta(u)'' 
\eeq
implies that 
\beq
    \frac{\rd\vartheta(u + \omega_1)/\rd\tau}
         {\vartheta(u + \omega_1)} 
    = \frac{1}{4\pi}
        \frac{\vartheta(u + \omega_1)''}
             {\vartheta(u + \omega_2)} 
    = \frac{1}{4\pi i} \left[ 
          \Bigl(\log\vartheta(u + \omega_1)\Bigr)'' 
          + \left(\frac{\vartheta'(u + \omega_1)}
            {\vartheta(u + \omega_1)} \right)^2 
        \right]. 
    \nonumber 
\eeq
By the third and forth lemmas, the last line 
can be rewritten 
\beq
    \frac{1}{4\pi i} \left[ - \wp(u + \omega_3) 
      + \left(2\pi i\frac{f_\tau(u)}{f'(u)}\right)^2 
    \right] 
    + \mbox{\rm function of $\tau$ only}, 
    \nonumber
\eeq
so that 
\beq
    2\pi i\frac{f_\tau(u)}{f'(u)}du 
    = \frac{1}{4\pi i}\left[ \wp(u + \omega_3) 
        - \left(2\pi i \frac{f_\tau(u)}{f'(u)}\right)^2 
      \right] d\tau 
    + \mbox{\rm exact form}. 
\eeq
Substituting $u = q$ gives (\ref{eq:term-A})

\section{Asymptotics of elliptic functions} 
\setcounter{equation}{0}

The asymptotic behavior of the $\wp$-function 
$\wp(u)$, the shifted $\wp$-functions 
$\wp(u + \omega_k)$ and the constants 
$e_k = \wp(\omega_k)$, in the limit as 
$\Im\tau \to +\infty$, can be deduced from 
the well known formula 
\beq
    \wp(u) 
    = \sum_{n=-\infty}^\infty 
      \frac{\pi^2}{\sin^2(\pi(u + n\tau))} 
    - \frac{\pi^2}{3} 
    - \sum_{n=1}^\infty 
      \frac{2\pi^2}{\sin^2(\pi n\tau)}. 
\eeq

Let us first consider the asymptotic behavior 
of $\wp(u)$ itself.  The constant ($n = 0$) 
term in the first sum is of order $1$ and 
the $n$-th term is of order $e^{2n\pi i\tau}$.  
Similarly, the $n$-th term in the second sum 
is of order $e^{2n\pi i\tau}$.  Therefore 
\beq
    \wp(u) = \frac{\pi^2}{\sin^2(\pi u)} 
           - \frac{\pi^2}{3} 
           + O(e^{2\pi i\tau}). 
\eeq

A similar estimate leads to the following 
asymptotic expression for the shifted 
$\wp$-functions: 
\beq
    \wp(u + \omega_1) 
    &=& \frac{\pi^2}{\cos^2(\pi u)} 
      - \frac{\pi^2}{3} 
      + O(e^{2\pi i\tau}), 
    \nonumber \\
    \wp(u + \omega_2) 
    &=& - \frac{\pi^2}{3} 
      + 8\pi^2 \cos(2\pi u) e^{\pi i\tau} 
      + O(e^{2\pi i\tau}), 
    \nonumber \\
    \wp(u + \omega_3) 
    &=& - \frac{\pi^2}{3} 
      - 8\pi^2 \cos(2\pi u)e^{2\pi i\tau} 
      + O(e^{2\pi i\tau}). 
\eeq
In fact, the degeneration process of the 
elliptic model requires us to know the 
asymptotic expression of $\wp(u + \omega_2) 
+ \wp(u + \omega_3)$ to the order 
$e^{2\pi i\tau}$.  This can be achieved 
by the following calculations: 
\beq
    && \wp(u + \omega_2) + \wp(u + \omega_3) 
    \nonumber \\
    &=& \sum_{n=-\infty}^\infty 
        \frac{\pi^2}
        {\cos^2(u + \frac{\tau}{2} + n\tau) 
         \sin^2(u + \frac{\tau}{2} + n\tau)} 
      - \frac{2\pi^2}{3} 
      - \sum_{n=1}^\infty 
        \frac{4\pi^2}{\sin^2(\pi n\tau)} 
    \nonumber \\
    &=& - \frac{2\pi^2}{3} 
       - 32\pi^2\cos(2\pi u) e^{2\pi i\tau} 
       + 16\pi^2 e^{2\pi i\tau} 
       + O(e^{3\pi i\tau}). 
\eeq

We now consider the constants $e_k$. 
For instance, $e_1$ can be written 
\beq
    e_1 
    &=& \sum_{n=-\infty}^\infty 
        \frac{\pi^2}{\cos^2(\pi n\tau)} 
      - \frac{\pi^2}{3} 
      - \sum_{n=1}^\infty 
        \frac{2\pi^2}{\sin^2(\pi n\tau)} 
    \nonumber \\
    &=& \frac{2}{3}\pi^2 
      + \sum_{n=1}^\infty 
        \frac{2\pi^2}{\cos^2(\pi n\tau)} 
      - \sum_{n=1}^\infty 
        \frac{2\pi^2}{\sin^2(\pi n\tau)}. 
\eeq
The constant $2\pi^2/3$ becomes the leading 
term; the leading ($n = 1$) terms of the last 
two series give the next-leading term of the 
order $e^{2\pi i\tau}$.  $e_2$ and $e_3$ can 
be similarly analyzed.  Thus the following 
asymptotic formulas are obtained: 
\beq
    e_1 &=& \frac{2\pi^2}{3} 
          + 16\pi^2 e^{2\pi i\tau} 
          + O(e^{4\pi i\tau}), 
    \nonumber \\
    e_2 &=& - \frac{\pi^2}{3} 
          + 8\pi^2 e^{\pi i\tau} 
          + O(e^{2\pi i\tau}) 
    \nonumber \\
    e_3 &=& - \frac{\pi^2}{3} 
          - 8\pi^2 e^{\pi i\tau} 
          + O(e^{2\pi i\tau}). 
\eeq
In particular, $e_2 - e_1 \to - \pi^2$, 
as expected.

%%%%%%%%%%%%%%%%%%%%
%%%% References %%%%
%%%%%%%%%%%%%%%%%%%%

\end{document}